\documentclass[11pt]{amsart}

\usepackage{amscd,amsmath,amssymb}
\usepackage[mathscr]{eucal}
\usepackage[all]{xy}

\makeatletter \@addtoreset{equation}{section}\makeatother

\newtheorem{theorem}{Theorem}[section]
\newtheorem{lemma}[theorem]{Lemma}
\newtheorem{proposition}[theorem]{Proposition}
\newtheorem{corollary}[theorem]{Corollary}

\newcommand{\op}{{\mathsf{op}}}
\newcommand{\sh}{{\mathsf{sh}}}
\newcommand{\eu}{{\mathsf{Eu}}}
\newcommand{\ch}{{\mathsf{Ch}}}
\newcommand{\str}{{\mathsf{str}}}
\newcommand{\tr}{{\mathsf{tr}}}
\newcommand{\eul}{{\mathsf{eu}}}

\newcommand{\Z}{{\mathbb{Z}}}
\newcommand{\C}{{\mathbb{C}}}
\newcommand{\X}{{\mathbb{X}}}

\newcommand{\Mod}{{\mathsf{Mod}}}
\newcommand{\free}{{\mathsf{Free}}}
\newcommand{\cone}{{\mathsf{Cone}}}
\newcommand{\Ho}{{\mathsf{Ho}}}

\newcommand{\Tw}{{\mathsf{Tw}}}

\newcommand{\K}{{\mathsf{K}}}
\newcommand{\per}{{\mathsf{Perf}}}

\newcommand{\perf}{{\mathsf{perf}}}

\newcommand{\qcoh}{{\mathsf{qcoh}}}
\newcommand{\mmod}{{\mathsf{mod}}}
\newcommand{\Hom}{{\mathrm{Hom}}}
\newcommand{\Harm}{{\mathrm{Harm}}}
\newcommand{\End}{{\mathrm{End}}}
\newcommand{\Ext}{{\mathrm{Ext}}}
\newcommand{\rHH}{{\mathsf{HH}}}
\newcommand{\rC}{{\mathsf{C}}}

\newcommand{\bHom}{{\mathbf{Hom}}}
\newcommand{\sFun}{{\mathsf{Fun}}}

\newcommand{\rH}{{\mathrm{H}}}
\newcommand{\bH}{{\mathbf{H}}}

\newcommand{\cA}{{\mathcal{A}}}
\newcommand{\cB}{{\mathcal{B}}}

\newcommand{\cF}{{\mathcal{F}}}
\newcommand{\cM}{{\mathcal{M}}}
\newcommand{\cV}{{\mathcal{V}}}

\title{\bf Hirzebruch-Riemann-Roch theorem for DG algebras}
\author{D. Shklyarov}

\begin{document}

\maketitle

\begin{center}{\it Dedicated to the memory of L. L. Vaksman}\end{center}

\baselineskip 1.2pc

\bigskip

\section{Introduction}
\subsection{Geometry of DG categories}\label{nagvdgc}
To motivate the subject of the present research, we will begin by
discussing some applications of triangulated and differential graded
categories in algebraic geometry.

Let $X$ be a quasi-compact separated scheme.\footnote{In what
follows, everything is considered over a fixed ground field.}
 Denote by $D_{\qcoh}(X)$ the derived category of complexes of
$\mathcal{O}_X$-modules with quasi-coherent cohomology and by
$D_{\perf}(X)$ its triangulated subcategory of perfect complexes,
i.e. complexes which are locally quasi-isomorphic to finite
complexes of vector bundles. The category $D_{\perf}(X)$ has proved
to be the basic (co)homological invariant of $X$ which somehow
encodes all other reasonable invariants. This idea underlies R.
Thomason's research on the $K$-theory of schemes \cite{TT}, M.
Kontsevich's Homological Mirror Symmetry program \cite{Kon0}, and A.
Bondal's and D. Orlov's research on the derived categories of smooth
schemes \cite{BO}.\footnote{One of their results claims that schemes
of certain type can be completely reconstructed from their derived
categories.}

When working with $D_{\perf}(X)$, one faces the following problem:
even though various invariants of $X$ depend on this category, it is
not clear how to compute some of them in terms of $D_{\perf}(X)$,
viewed as an abstract triangulated category. One way to get around
the problem is due to A. Bondal and M. Kapranov \cite{BK}. The point
is that the derived categories, as opposed to abstract triangulated
categories, can be ``upgraded" to differential graded (DG)
categories. In practice this can be achieved by, say, passing from
$D_{\perf}(X)$ to the DG category $\per X$ of left bounded injective
perfect complexes. The category $D_{\perf}(X)$ is then recovered as
the homotopy category of $\per X$. Many other invariants of $X$ can
be extracted from $\per X$ as well. The simplest example is the
computation of the Hodge cohomology of $X$ in terms of $\per X$ in
the case of a smooth scheme. One has
\begin{equation}\label{hodge}
\rHH_n(\per X)=\oplus_i \rH^{i-n}(\Omega^i_X),
\end{equation}
where the left-hand side stands for the $n$-th Hochschild homology
group of $\per X$ (see Section \ref{hoch}). The category $\per X$
encodes also some geometric properties of $X$. For example, if $X$
is smooth then the category $\per X$ is a perfect bimodule over
itself \cite{KS}.

The DG categories of the form $\per X$ turn out to be equivalent to
the DG categories of perfect modules over certain DG algebras.
Namely, according to \cite[Section 3.1]{BVDB} (see also \cite{Rouq})
$D_{\perf}(X)$ is generated by a single perfect complex,
$\mathcal{E}$. Let $A=\End_{\per X}(\mathcal{E})$. Then $\per X$ is
quasi-equivalent to the DG category $\per A$ of perfect right
$A$-modules (see Section \ref{perf} for the definition of the latter
category). Of course, there is no canonical generator of
$D_{\perf}(X)$ and, as a result, there is no canonical DG algebra
associated with the scheme. However any DG algebra such that $\per
X$ is quasi-equivalent to $\per A$ can be viewed as a replacement of
the algebra of regular functions in the case of a non-affine scheme
$X$.

Let us look at the most popular example - the projective line
$\mathbf{P}^1$. Due to the well known result of A. Beilinson
\cite{Be}, the derived category of coherent sheaves in this case is
equivalent to the derived category of finite dimensional modules
over the path algebra of the Kronecker quiver:
\begin{displaymath}
\xymatrix{ \bullet \ar@/^/[r] \ar@/_/[r] & \bullet }
\end{displaymath}

Following \cite{To}, we will say that two DG algebras $A$ and $B$
are Morita-equivalent if their perfect categories $\per A$ and $\per
B$ are quasi-equivalent. In view of the above discussion, each
scheme gives rise to a fixed Morita-equivalence class. Therefore it
is reasonable to think of an {\it arbitrary} Morita-equivalence
class as representing some noncommutative scheme or, better yet, a
{\it noncommutative DG-scheme}. Any DG algebra from the equivalence
class should be viewed as ``the" algebra of regular functions on
this noncommutative DG-scheme, and $\per A$ plays the role of $\per
X$.

The above point of view agrees with the philosophy of {\it derived
noncommutative algebraic geometry}.\footnote{We are not sure whether
this name is commonly accepted or not.} This subject was initiated
in the beginning of 90's based on the previous extensive study of
derived categories of coherent sheaves undertaken by the Moscow
school (A. Beilinson, A. Bondal, M. Kapranov, D. Orlov, A. Rudakov
et al). Later on, it was greatly enriched by new ideas and examples
coming from M. Kontsevich's Homological Mirror Symmetry program
\cite{Ko}. A particularly important implication of the program is
that one can associate certain triangulated categories with
symplectic manifolds which should play the same important role in
symplectic geometry that the derived categories of coherent sheaves
play in algebraic geometry. Further important ideas and results in
the field are due to A. Bondal and M. Van den Bergh, T. Bridgeland,
V. Drinfeld, B. Keller, M. Kontsevich and Y. Soibelman, D. Orlov, R.
Rouquier, B. Toen and others.

Of course, a ``real" definition of noncommutative DG-schemes should
include also a description of morphisms between them. It is clear
that morphisms are given by DG functors between the categories of
perfect complexes (a prototype is the pull-back functor associated
with a morphism of schemes). The real definition is more subtle and
we won't discuss it here referring the reader to more thorough
treatments of the subject \cite{D,K2,Tab,Ta,To,TV}.

Here is an interesting question: Is it possible to tell whether a
noncommutative DG-scheme comes from a usual commutative one? There
is a simple necessary condition: the corresponding DG algebra $A$
should be Morita-equivalent to its opposite DG algebra $A^\op$  (the
simplest case when this is so is when the DG algebras $A$ and
$A^\op$ are isomorphic; look at the Kronecker quiver!). Of course,
this condition is not sufficient: various almost commutative
schemes, such as orbifolds, also satisfy it.

Let $A$ be a DG algebra. Then, following \cite{KS}, one can define
the corresponding noncommutative DG-scheme to be

{\it proper} iff so is $A$, i.e. $\sum_n\dim\rH^n(A)<\infty$;

{\it smooth} iff so is $A$, i.e. $A$ is quasi-isomorphic to a
perfect $A$-bimodule.

\noindent The first property is central to the present work,
although we will touch upon smooth DG algebras as well (see Section
\ref{pairi}).

\medskip
\subsection{A categorical version of the Hirzebruch-Riemann-Roch theorem}\label{chrr}
Let us turn now to the subject of this article,
Hirzebruch-Riemann-Roch (HRR) theorem in the above noncommutative
setting. We will start with very general (and oversimplified)
categorical considerations.

Fix a ground field, $k$, and consider the tensor category of small
$k$-linear DG categories, morphisms being DG functors. Fix also a
{\it homology theory} on the latter category, i.e. a covariant
tensor functor $\bH$ to a tensor category of modules over a
commutative ring\footnote{One can take $\Z$-graded, $\Z/2$-graded
modules, modules that are complete in some topology etc.} $K$,
satisfying the following axioms:

\medskip

(1) $\bH$ respects quasi-equivalences.

(2) For any DG algebra $A$ the canonical embedding $A\to\per A$
induces an isomorphism $$\bH(A)\simeq\bH(\per A).$$

(3) $\bH(k)=K$ (then, by (2),\, $\bH(\per k)=K$).

\medskip

Notice that (1) and (2) together imply that $\bH$ descends to an
invariant of noncommutative DG-schemes. Also, by the very definition
of $\bH$, there exists a functorial K\"{u}nneth type isomorphism
$$\bH(\cA)\otimes_K\bH(\cB)\simeq\bH(\cA\otimes\cB).$$

Let us add to this list one more condition:

\medskip

(4) For any DG category $\cA$ there is a functorial isomorphism
$$\,^\vee: \bH(\cA)\simeq\bH(\cA^\op)$$ which equals identity in the
case $\cA=k$.

\medskip

We will assume that the above isomorphisms satisfy all the natural
properties and compatibility conditions one can imagine
\footnote{The right definition of a homology theory should be
formulated in terms of the category of noncommutative motives
\cite{Kon}.}.

To describe what we understand by an abstract HRR theorem for
noncommutative DG-schemes, we need to define the Chern character map
with values in the homology theory $\bH$. This is a function
$\ch^\cA_\bH: \cA\to\bH(\cA)$, one for each DG category $\cA$,
defined as follows. Take an object $N\in\cA$ and consider the DG
functor $T_N: k\to\cA$ that sends the unique object of $k$ to $N$.
Then \cite{BNT,K3}
$$
\ch^\cA_\bH(N)=\bH(T_N)(1_K).
$$
Clearly, the Chern character is functorial: For any two DG
categories $\cA,\cB$ and any DG functor $F:\cA\to\cB$
$$
\ch^\cB_\bH\circ F=\bH(F)\circ\ch^\cA_\bH.
$$

From now on, we will focus on proper DG categories, i.e. DG
categories that correspond to proper noncommutative DG-schemes. Let
$\cA$ be a proper DG category. Consider the DG functor
$$
\bHom_\cA: \cA\otimes\cA^\op\to\per k,\quad N\otimes M\mapsto\Hom_{\cA}(M,N).
$$
By (3), it induces a linear map $\bH(\bHom_\cA):
\bH(\cA\otimes\cA^\op)\to K$. One can compose it with the
K\"{u}nneth isomorphism to get a $K$-bilinear pairing
$$
\langle\,,\,\rangle_\cA: \bH(\cA)\times\bH(\cA^\op)\to K.
$$

Now we are ready to formulate the HRR theorem: For any proper DG
category $\cA$ and any two objects $N,M\in\cA$
\begin{equation}\label{crr}
\ch^{\per k}_\bH(\Hom_{\cA}(M,N))=\langle\,\ch^\cA_\bH(N)\,,\,\ch^\cA_\bH(M)^\vee\,\rangle_\cA.
\end{equation}
Indeed, it follows from the functoriality of the isomorphism
$\,^\vee$ that $$(\bH(T_M)(1_K))^\vee=\bH(T_{M^\op})(1_K)$$ where
$M^\op$ stands for $M$ viewed as an object of $\cA^\op$. Then
\begin{eqnarray*}
\langle\,\ch^\cA_\bH(N)\,,\,\ch^\cA_\bH(M)^\vee\,\rangle_\cA=
\bH(\bHom_\cA) \left(\bH(T_N)(1_K)
\otimes(\bH(T_M)(1_K))^\vee\right)\\
=\bH(\bHom_\cA)\left(\bH(T_N)(1_K)\otimes\bH(T_{M^\op})(1_K)\right)=\bH(\bHom_\cA)\left(\bH(T_{N\otimes
M^\op})(1_K)\right)\\=\bH(\bHom_\cA\circ T_{N\otimes
M^\op})(1_K)=\bH(\Hom_{\cA}(M,N))(1_K)=\ch^{\per
k}_\bH(\Hom_{\cA}(M,N)).
\end{eqnarray*}

In this very general form, the HRR theorem is almost tautological.
For it to be of any use, one needs to find a way to compute the
right-hand side of (\ref{crr}) for a given proper noncommutative
DG-scheme and any pair of perfect complexes on it. In this work, we
solve this problem in the case $K=k$, $\bH=\rHH_\bullet$, where
$\rHH_\bullet$ stands for the Hochschild homology\footnote{The most
difficult axioms (1) and (2) in our ``definition" of the homology
theory were proved for $\rHH_\bullet$ by B. Keller in \cite{K1}.}
(see Section \ref{hoch} for the definition). This choice of the
homology theory can be motivated as follows.

First of all, there is a classical character map from the
Grothendieck group of a ring to its Hochschild homology - the so
called Dennis trace map \cite{L}. Its sheafified version appeared in
\cite{BNT} in connection with the index theorem for elliptic pairs
\cite{SS1,SS2} (the definition of the Chern character given above
mimics the one given in \cite{BNT}).

In the algebraic geometric context, the relevance of the Hochschild
homology to the HRR theorem can be explained as follows. There is a
version of the HRR theorem for compact complex manifolds
\cite{OBTT1,OBTT2}, in which the Chern class of a coherent sheaf
takes values in the Hodge cohomology $\oplus_i \rH^i(\Omega^i_X)$
(see also \cite{I}). A new proof of this result was obtained in
\cite{Ma1,Ma2} using an algebraic-differential calculus (see also
\cite{Ca1,Ra}). This latter approach emphasizes the importance of
viewing the Chern character as a map to the Hochschild homology
$\rHH_0(X)$ of the space $X$. The ``usual" Chern character is then
obtained via the Hochschild-Kostant-Rosenberg isomorphism
$\rHH_0(X)\cong\oplus_i\rH^i(\Omega^i_X)$. This point of view was
further developed in \cite{Ca}. Namely, it was explained in
\cite{Ca} (see also \cite{CaWi}) how to obtain a categorical version
of the HRR theorem, similar to the one above, starting from the
cohomology theory
$$
\text{smooth}\,\text{spaces}\to\text{graded}\,\text{vector}\,\text{spaces},\quad X\mapsto\rHH_\bullet(X)
$$
(``smooth spaces" are understood in a broad sense: these are usual
schemes as well as various almost commutative ones such as
orbifolds). Finally, the transition from $X$ to its categorical
incarnation, $\per X$, is based on the fact that $\rHH_\bullet(X)$
is isomorphic to the Hochschild homology of the DG category $\per
X$, which was proved in \cite{K3}.

Before we move on to the description of the main results of the
paper, we would like to mention a notational convention we are going
to follow.

Following \cite{BNT} (see also \cite{K3}), we will call the Chern
character $\ch_\rHH$ with values in the Hochschild homology the {\it
Euler} character and use the notation $\eu$.

\medskip
\subsection{Main results}
Let us describe the main results of this work.

Fix a ground field $k$ and a proper DG algebra $A$ over $k$ (as we
mentioned earlier, the properness means
$\sum_n\dim\rH^n(A)<\infty$).

The first main result is the computation of the Euler class
$\eul(L)$ of an arbitrary perfect DG $A$-module $L$. Here $\eul(L)$
stands for the unique element in $\rHH_0(A)$ that corresponds to
$\eu(L)\in\rHH_0(\per A)$ under the canonical isomorphism
$\rHH_\bullet(A)\simeq\rHH_\bullet(\per A)$ (see axiom (2) in
Section \ref{chrr}). The following theorem is proved in Section
\ref{cec}.

\medskip

\begin{quote}
{\bf Theorem 1.} \emph{Let $N=(\bigoplus_jA[r_j], d+\alpha)$ be a
twisted DG $A$-module and $L$ a homotopy direct summand of $N$ which
corresponds to a homotopy idempotent $\pi:N\to N$. Then
$$
\eul(L)=\sum_{l=0}^\infty
(-1)^l\str(\pi[\underbrace{\alpha|\ldots|\alpha}_{l}])
$$
}
\end{quote}

Roughly speaking, in this formula $\pi$ and $\alpha$ are elements of
a DG analog of the matrix algebra $\mathsf{Mat}(A)$,
$\pi[\alpha|\ldots|\alpha]$ is an element of the Hochschild chain
complex of this DG matrix algebra, and $\str$ is an analog of the
usual trace map $\tr:\mathsf{Mat}(A)\to A$ (see \cite{Gi,L}). Note
that $\alpha$ is upper-triangular, so the series terminates.

\medskip

To present our next result, we observe that the pairing
$$
\rHH_\bullet(\per A)\times\rHH_\bullet((\per A)^\op)\to \rHH_\bullet(\per k)\simeq k,
$$
defined earlier in Section \ref{chrr}, induces a pairing
\begin{equation}\label{pp}
\rHH_\bullet(\per A)\times\rHH_\bullet(\per A^\op)\to k.
\end{equation}
This is due to the existence of a canonical quasi-equivalence of DG
categories (see (\ref{qe})):
$$
D: \per A^\op\to(\per A)^\op, \quad M\mapsto DM=\Hom_{\per A^\op}(M, A).
$$
In fact, we ``twist" the exposition in the main text (Section
\ref{ec}) and work exclusively with the pairing (\ref{pp}). The
reason is that it can be defined very explicitly without referring
to its categorical nature. Besides, it induces a pairing
\begin{equation}\label{hhp}
\langle\,,\,\rangle: \rHH_\bullet(A)\times\rHH_\bullet(A^\op)\to k
\end{equation}
via the canonical isomorphisms
$\rHH_\bullet(A)\simeq\rHH_\bullet(\per A)$,
$\rHH_\bullet(A^\op)\simeq\rHH_\bullet(\per A^\op)$. This latter
pairing is described explicitly in our next theorem, which is
obtained by combining results of Section \ref{hrr} (see formulas
(\ref{wedge}), (\ref{pairing})) and Theorem \ref{integral}.

\medskip

\begin{quote}
{\bf Theorem 2.} \emph{Let $a$, $b$ be two elements of
$\rHH_\bullet(A)$, $\rHH_\bullet(A^\op)$, respectively. Then
$$
\langle\,a,b\,\rangle=\int a\wedge b.
$$
Here $\wedge: \rHH_\bullet(A)\times\rHH_\bullet(A^\op)\to
\rHH_\bullet(\End_k(A))$, $\int:\rHH_\bullet(\End_k(A))\to k$ are
defined as follows:}

\medskip

\emph{\noindent(1) If $\sum_{a}a_0[a_1|\ldots |a_l]$ (resp.
$\sum_{b}b_0[b_1|\ldots |b_m]$) is a cycle in the Hochschild chain
complex of $A$ (resp. $A^\op$) representing the homology class $a$
(resp. $b$) then
$$
a\wedge b=\sum_{a,b}\sh\left(L(a_0)[L(a_1)|\ldots |L(a_l)]\otimes R(b_0)[R(b_1)|\ldots |R(b_m)]\right),
$$
where $L(a_i)$ (resp. $R(b_j)$) stands for the operator in $A$ of
left (resp. right) multiplication with $a_i$ (resp. $b_j$); $\sh$ is
the well known shuffle-product (see Section \ref{Ku}).}

\medskip

\emph{\noindent (2) $\int$ is what we call the Feigin-Losev-Shoikhet
trace \cite{FLS,R}. It is described explicitly in Theorem
\ref{integral} (Section \ref{ci}).}
\end{quote}

\medskip

Furthermore, recall that there should exist a canonical isomorphism
$\,^\vee: \rHH_\bullet(A)\simeq\rHH_\bullet(A^\op)$ (see axiom (4)
in Section \ref{chrr}). In fact, the isomorphism is easy to describe
explicitly (see Section \ref{hrr}). By summarizing the above
discussion, we get the following version of the noncommutative HRR
theorem:

\medskip

\begin{quote}
{\bf Theorem 3.} \emph{For any perfect DG $A$-modules $N,M$
$$\chi(M,N)(:=\chi(\Hom_{\per A}(M,N)))=\int\eul(N)\wedge\eul(M)^\vee.$$
}
\end{quote}
The only thing that needs to be explained here is where
$\chi(\Hom_{\per A}(M,N))$ came from. According the  categorical HRR
theorem, described in the previous section, the left-hand of the
above equality should equal $\eul(\Hom_{\per A}(M,N))$. However, the
Euler class of a perfect DG $k$-module is nothing but its Euler
characteristic. This is a consequence of the following ``expected"
fact, which we prove in Section \ref{ec}: for any $A$ the Euler
character $\eu$ descends to a character on the Grothendieck group of
the triangulated category $\Ho(\per A)$, the homotopy category of
$\per A$.

\medskip

Note that the noncommutative HRR formula doesn't include any analog
of the Todd class. The Todd class seems to emerge in the case when a
noncommutative space, $\widehat{X}$, is ``close" to a commutative
one, $X$ (for example, $\widehat{X}$ is a deformation quantization
of $X$). In such cases various homology theories of $\widehat{X}$
can be identified with certain cohomology rings associated with $X$
and the Todd class of $X$ appears because of this identification.
For some classes of noncommutative spaces one can try to define an
analog of the Todd class ``by hand"  but, in general, a categorical
definition of the Todd class doesn't seem to exist.

In the main text we do not refer to the categorical version of the
HRR theorem to prove Theorem 3. Instead, we derive it from a more
general statement (Theorem \ref{main}). Roughly speaking, this
statement says the following: If $A$ and $B$ are two proper DG
algebras and $X$ is a perfect $A-B$-bimodule then the map
$\rHH_\bullet(\per A)\to\rHH_\bullet(\per B)$, induced by the DG
functor $-\otimes_AX:\per A\to\per B$, is given by a ``convolution"
with the Euler class of $X$. Later, in Section \ref{smooth}, we use
this result again to prove the following

\medskip

\begin{quote}
{\bf Theorem 4.} \emph{Let $A$ be a proper smooth DG algebra. Then
the pairing $\langle\,,\,\rangle$ is non-degenerate. }
\end{quote}
It is this application that was the original motivation for the
author to study the Euler classes in the DG setting \footnote{I am
grateful to Y. Soibelman for suggesting to me to ``write up" the
proof of this statement.}. It implies, in particular, the
noncommutative Hodge-to-De Rham degeneration conjecture \cite{KS}
for smooth algebras with the trivial differential and grading.
Hopefully, the reader will accept all this as an excuse for
``twisting" the exposition and not mentioning the categorical HRR in
what follows.

\medskip

In the end of this work we present some ``toy" examples of proper
noncommutative DG-schemes and the HRR formulas for them. Namely, in
Section \ref{dir} we discuss what we call directed algebras.
Basically, these are some quiver algebras with relations but we find
the quiver-free description more convenient when it comes to proving
general facts about such algebras. Many commutative schemes, viewed
as noncommutative ones, are described by directed algebras. Namely,
this is so when the scheme possesses a strongly exceptional
collection \cite{B}. The HRR formula for such algebras (see
(\ref{grf})) is a special case of Ringel's formula \cite[Section
2.4]{Ri}. Section \ref{orbi} is about proper noncommutative
DG-schemes ``responsible" for orbifold singularities of the form
$\C^n/G$, where $G$ is a finite subgroup of $SL_n(\C)$. Namely, we
look at the noncommutative DG-scheme related to the derived category
of complexes of $G$-equivariant coherent sheaves on $\C^n$ with
supports at the origin. We conjecture that the underlying DG algebra
is the cross-product $\Lambda^\bullet\C^n\rtimes\C[G]$ and we derive
the HRR formula for some perfect modules over this algebra (see
(\ref{orbhrr})).

Section \ref{TCFT} is devoted to a less straightforward application
of our results. It has been conjectured by Y. Soibelman and K.
Costello that for a Calabi-Yau DG algebra the pairing
$\rHH_\bullet(A)\times \rHH_\bullet(A^\op)\to k$, we construct in
this paper, coincides with the one coming from the Topological Field
Theory associated with $A$ by \cite{Cos,KS}. In Section \ref{TCFT}
we formulate this conjecture and verify it in the particular case of
Frobenius algebras.

\medskip
\subsection{Other viewpoints on the noncommutative HRR theorem}
In this section, we provide a very brief account of other
Riemann-Roch type theorems in Noncommutative Geometry we are aware
of.

Let us begin with the afore-mentioned preprint \cite{Ca} which
partially inspired the present work. The approach taken in \cite{Ca}
is based on an alternative description of the Hochschild homology of
a {\it smooth} proper space $X$ in terms of the Serre functor $S_X:
D^b(X)\to D^b(X)$. Namely,
$$
\rHH_\bullet(X)\simeq\Ext^\bullet_{\sFun}(S_X^{-1}, I_X),
$$
where $I_X$ is the identity endofunctor of $D^b(X)$ and the
extensions are taken in a suitably defined triangulated category of
endofunctors. In \cite{Sh} we generalized the above isomorphism to
the case of an arbitrary smooth proper noncommutative DG-scheme.
However, proving that the above definition gives rise to a homology
theory on the category of smooth proper noncommutative DG-schemes
(in other words, lifting the above definition on the level of DG
categories) will require some efforts \cite[Appendix B]{Ca}.
Besides, the ``traditional" definition of the Hochschild homology we
use in this paper works for an arbitrary, not necessarily smooth
scheme.

Other analogs of the Riemann-Roch theorem were obtained in \cite{J},
\cite{Mo} in the framework of Noncommutative Algebraic Geometry
\cite{AZ}, \cite{Ros1}, \cite{Ros2}, \cite{VdB1}. The exposition in
\cite{Mo} is closer to ours in that it emphasizes the importance of
triangulated categories in connection with Riemann-Roch type
results. Our approach and the above two approaches to the
noncommutative Riemann-Roch theorem are not completely unrelated
since many interesting noncommutative schemes give rise to
noncommutative DG schemes \cite[Section 4]{BVDB}.

Last, but not least, various index theorems have been proved in
frameworks of A. Connes' Noncommutative Geometry \cite{Con,Sk} and
Deformation Quantization \cite{BNT,F,NT,Ts}.

\medskip
\subsection{Acknowledgements} I am indebted to Y. Soibelman for prompting my
interest in derived noncommutative algebraic geometry and numerous
discussions which have played a crucial role in shaping my
understanding of the subject. I am also very grateful to B. Keller
and D. Orlov for patiently answering my questions about triangulated
and DG categories, to K. Costello for numerous interesting remarks
and encouragement, and to B. Tsygan for interesting comments on
noncommutative Chern characters and the noncommutative Hodge-to-De
Rham degeneration conjecture. Finally, I would also like to thank B.
Keller, Yu. I. Manin, A. L. Rosenberg, M. Van den Bergh for sending
me their comments which helped me to improve the exposition.
Needless to say, I am alone responsible for typos and more serious
mistakes if there are any.

This work was partially supported by NSF grant DMS-0504048.

\medskip
\section{Hochschild homology of DG algebras and DG categories}
\subsection{DG algebras, DG categories, and DG functors}\label{1}

Throughout the paper, we work over a fixed ground field $k$. All
vector spaces, algebras, linear categories are defined over $k$.

We consider unital DG algebras with no restrictions on the
$\Z$-grading. If $A$ is a DG algebra
$$
A=\bigoplus\limits_{n\in\mathbb{Z}}A^n,\quad d=d_A: A^n\rightarrow
A^{n+1}
$$
then $A^\op$ will stand for the opposite DG algebra, i.e. $A^\op$
coincides with $A$ as a complex of vector spaces and the product on
$A^\op$ is given by
$$
a'\otimes a''\mapsto(-1)^{|a'||a''|}a''a'
$$
(here and further, $|n|$ denotes the degree of a homogeneous element
$n$ of a graded vector space). $\Mod A$ will stand for the DG
category of right DG $A$-modules.

The homotopy category of a DG category $\cA$ will be denoted by
$\Ho(\cA)$. Let us recall the definition of the standard
triangulated structure on $\Ho(\Mod A)$. The shift functor is
defined in the obvious way. The distinguished triangles are defined
as follows. Let $p: L\to M$ be a degree 0 closed morphism. The cone
$\cone(p)$ of the morphism $p$ is a DG $A$-module defined by
$$
\cone(p)=\left(\begin{array}{c}L[1]\\ \oplus\\ M\end{array}, \begin{pmatrix}d_{L[1]}
& 0 \\ p & d_M \end{pmatrix}\right)
$$
(the direct sum is taken in the category of {\it graded}
$A$-modules). There are obvious degree 0 closed morphisms $q: M\to
\cone(p)$ and $r: \cone(p)\to L[1]$. A triangle in $\Ho(\Mod A)$ is,
by definition, a sequence $X\to Y\to Z\to X[1]$ isomorphic to a
sequence of the form $L\stackrel{p}\to M\stackrel{q}\to \cone(p)
\stackrel{r}\to L[1]$.

Let $N$ be a right DG $A$-module. A degree 0 closed endomorphism
$\pi: N\to N$ will be called a homotopy idempotent if $\pi^2=\pi$ in
$\Ho(\Mod A)$. By a homotopy direct summand of $N$ we will
understand a DG $A$-module $L$ that satisfies the following
property: there exists a homotopy idempotent $\pi: N\to N$ and two
degree 0 morphisms $f:N\to L$ and $g:L\to N$ such that $fg=1_{L}$,
$gf=\pi$ in $\Ho(\Mod A)$.

Fix two DG categories $\cA$ and $\cB$ and consider the DG category
$\sFun(\cA,\cB)$ of DG functors from $\cA$ to $\cB$ \cite{K2}. A
degree 0 closed morphism $f\in\Hom_{\sFun(\cA,\cB)}(F, G)$ will be
called a {\it weak homotopy equivalence} if for any $N\in\cA$ the
morphism $f(N):F(N)\to G(N)$ is a homotopy equivalence in
$\Ho(\cB)$.

\medskip
\subsection{Perfect modules}\label{perf}
Let $A$ be a DG algebra. It can be viewed as a full DG subcategory
of $\Mod A$ with a single object. The embedding
$A\hookrightarrow\Mod A$ factors through a smaller full subcategory
$\per A\subset\Mod A$ of perfect $A$-modules. This subcategory is
defined as follows (see \cite{BLL}).

Let us say that a DG $A$-module $N$ is {\it finitely generated free}
if it is isomorphic to a module of the form $K\otimes A$ where $K$
is a finite dimensional graded vector space (equivalently, it is a
finite direct sum of shifts of $A$). We will say that $N\in\Mod A$
is {\it finitely generated semi-free} if it can be obtained from a
finite set of finitely generated free $A$-modules (equivalently, a
finite set of shifts of $A$) by successive taking the cones of
degree 0 closed morphisms. Finally, a {\it perfect} DG $A$-module is
a homotopy direct summand of a finitely generated semi-free DG
$A$-module.

Note that this definition is slightly more general than the one
given in \cite{BLL}. The authors of \cite{BLL} require perfect
modules to be semi-free but we don't. For example, a complex of
vector spaces is perfect in our sense iff it has finite dimensional
total cohomology and it is perfect in the sense of \cite{BLL} if, in
addition, it is bounded above. The reason we prefer not to restrict
ourselves to semi-free modules will be clear from Proposition
\ref{tensor} below. It suffices for us to stay within the class of
homotopically projective modules: $N$ is homotopically projective
iff $\Hom_{\Mod A}(N,L)$ is acyclic whenever $L$ is acyclic. Every
finitely generated semi-free module $N$ is known to be homotopically
projective \cite[Section 13]{D}. It follows that every perfect
module in our sense is homotopically projective as well.

The following result is well known (and is not hard to prove):

\medskip
\begin{proposition}\label{summ}
The DG category $\per A$ is closed under passing to homotopically
equivalent modules, taking shifts and cones of degree 0 morphisms,
and taking homotopy direct summands.
\end{proposition}

\medskip

Let us list some simple useful facts about DG functors between the
categories of perfect modules.

\medskip
\begin{proposition}\label{preper}
Let $A$, $B$ be DG algebras and $F:\Mod A\to\Mod B$ a DG functor.
The DG functor $F$ preserves the subcategories of perfect modules
iff $F(A)\in\per B$.
\end{proposition}
To prove this proposition, observe that $F$ preserves homotopy
direct summands and cones of degree 0 morphisms.

\medskip
For two DG algebras $A,B$ and a bimodule $X\in\Mod(A^\op\otimes B)$
let us denote by $T_X$ the DG functor $$-\otimes_A X:\Mod A\to\Mod
B.$$ Here is a straightforward consequence of the last proposition:

\medskip
\begin{corollary}\label{cor1}
Suppose a bimodule $X\in\Mod(A^\op\otimes B)$ is perfect as a DG
$B$-module. Then $T_X$ preserves perfect modules.
\end{corollary}

\medskip

Recall that a DG algebra $A$ is called {\it proper} if
$\sum_n\dim\rH^n(A)<\infty$.

\medskip

\begin{proposition}\label{tensor}
Let $A$ be a proper DG algebra and $B$ an arbitrary DG algebra. Then
for any $X\in\per(A^\op\otimes B)$ the DG functor $T_X$ preserves
perfect modules.
\end{proposition}
In view of the above corollary, it is enough to show that $X$ is
perfect as a DG $B$-module. Suppose that $X$ is a homotopy direct
summand of a finitely generated semi-free DG $A^\op\otimes B$-module
$Y$ and $Y$ is obtained from $(A^\op\otimes B)[m_1],\ldots,
(A^\op\otimes B)[m_l]$ by successive taking cones of degree 0 closed
morphisms. As a $B$-module, $A^\op\otimes B$ is homotopically
equivalent to the finitely generated free module
$\rH^\bullet(A)\otimes B$ (this is where we use the properness of
$A$ and the fact that we are working over a field!). Thus, as a
$B$-module, $Y$ is homotopically equivalent to a finitely generated
semi-free module. Then $X$, as a $B$-module, is a homotopy direct
summand of a module that is homotopy equivalent to a finitely
generated semi-free module. This, together with Proposition
\ref{summ}, finishes the proof.

\medskip

Let us recall one more result about perfect modules which we will
need later on. The fact that perfect modules are homotopically
projective implies the following result (cf. \cite[Corollary
10.12.4.4]{BL}):

\medskip

\begin{proposition}\label{flat}
If $P$ is a perfect right DG $A$-module then $P\otimes_AN$ is
acyclic for every acyclic DG $A^\op$-module $N$.
\end{proposition}

\medskip

\subsection{Hochschild homology}\label{hoch}

We begin by recalling the definition of the Hochschild homology
groups $\rHH_n(A)$, $n\in\Z$, of a DG algebra $A$.

Let us use the notation $sa$ to denote an element $a\in A$ viewed as
an element of the ``suspension" $sA=A[1]$. Thus, $|sa|=|a|-1$. Let
$\rC_\bullet(A)=A\otimes T(A[1])=\bigoplus\limits_{n=0}^\infty
A\otimes A[1]^{\otimes n}$ equipped with the induced grading. We
will denote elements of $A\otimes A[1]^{\otimes n}$ by $a_0$, if
$n=0$, and $a_0[a_1|a_2|\ldots |a_n]$ otherwise (i.e.
$a_0[a_1|a_2|\ldots |a_n]=a_0\otimes sa_1\otimes
sa_2\otimes\ldots\otimes sa_n$). $\rC_\bullet(A)$ is equipped with
the differential $b=b_0+b_1$, where $b_0$ and $b_1$ are two
anti-commuting differentials given by
\begin{equation}\label{d1d2}
b_0(a_0)=da_0,\quad b_1(a_0)=0,
\end{equation}
and
\begin{eqnarray*}\label{D1D2}
b_0(a_0[a_1|a_2|\ldots |a_n])=da_0[a_1|a_2|\ldots |a_n]
-\sum\limits_{i=1}^n(-1)^{\eta_{i-1}}a_0[a_1|a_2|\ldots |da_i|\ldots|a_n],
\end{eqnarray*}
\begin{eqnarray*}
b_1(a_0[a_1|a_2|\ldots |a_n])=(-1)^{|a_0|}a_0a_1[a_2|\ldots |a_n]
+\sum\limits_{i=1}^{n-1}(-1)^{\eta_{i}}a_0[a_1|a_2|\ldots
|a_ia_{i+1}|\ldots|a_n]\\
-(-1)^{\eta_{n-1}(|a_n|+1)}a_na_0[a_1|a_2|\ldots |a_{n-1}]
\end{eqnarray*}
for $n\neq0$.  Here $\eta_i=|a_0|+|sa_1|+\ldots+|sa_i|$.
$\rC_\bullet(A)$ is called the Hochschild chain complex of $A$. Then
$$
\rHH_n(A)=\rH^n(\rC_\bullet(A)).
$$

Let $\cA$ be a (small) DG category. Its Hochschild chain complex is
defined as follows. Fix a non-negative integer $n$. We will denote
the set of sequences $\{X_0,X_1,\ldots,X_n\}$ of objects of $\cA$ by
$\cA^{n+1}$ (the objects in the sequence are not required to be
different). Fix an element $\X=\{X_0,X_1,\ldots,X_n\}\in\cA^{n+1}$
and denote by $\rC_\bullet(\cA,\X)$ the graded vector space
$\Hom_{\cA}(X_n,X_0)\otimes
\Hom_{\cA}(X_{n-1},X_n)[1]\otimes\ldots\otimes
\Hom_{\cA}(X_0,X_1)[1]$. Equip the space
$$
\rC_\bullet(\cA)=\bigoplus_{n\geq0}\bigoplus_{\X\in\cA^{n+1}}\rC_\bullet(\cA,\X)
$$
with the differential $b=b_0+b_1$ where $b_0$ and $b_1$ are given by
formulas analogous to (\ref{d1d2}),(\ref{D1D2}). The complex
$\rC_\bullet(\cA)$ is the Hochschild chain complex of the DG
category $\cA$ and its cohomology
$$
\rHH_n(\cA)=\rH^n(\rC_\bullet(\cA))
$$
is the Hochschild homology of $\cA$.

Obviously, any DG functor $F:\cA\to\cB$ between two DG categories
$\cA$, $\cB$ induces a morphism of complexes
$\rC(F):\rC_\bullet(\cA)\to \rC_\bullet(\cB)$ and, as a result, a
linear map $$\rHH(F):\rHH_\bullet(\cA)\to \rHH_\bullet(\cB).$$

Being applied to the embedding $A\to\per A$, the above construction
yields a morphism of complexes $\rC_\bullet(A)\to \rC_\bullet(\per
A)$. The following result was proved in \cite{K1} (see also
\cite{K2}):

\medskip
\begin{theorem}\label{hhper}
The morphism $\rC_\bullet(A)\to \rC_\bullet(\per A)$ is a
quasi-isomorphism.
\end{theorem}
\medskip

Later on, we will need yet another result proved in \cite{K1} (see
Section 3.4 of loc.cit.):
\medskip
\begin{theorem}\label{hominv}
Let $A$ and $B$ be two DG algebras and $F,G:\per A\to\per B$ two DG
functors. If there is a weak homotopy equivalence $F\to G$ then
$\rHH(F)=\rHH(G)$.
\end{theorem}

\medskip
\subsection{K\"{u}nneth isomorphism}\label{Ku}

Let us recall the construction of the K\"{u}nneth isomorphism
$$\bigoplus_{n}\rHH_n(A)\otimes\rHH_{N-n}(B)\simeq\rHH_N(A\otimes B)$$ where
$A,B$ are two DG algebras. The formula below is borrowed from
\cite{L} (see also \cite{Ts} where the differential graded case is
discussed).

Let us fix a DG algebra $A$. The first ingredient of the
construction is the shuffle product
$$
\sh: \rC_\bullet(A)\otimes \rC_\bullet(A)\to \rC_\bullet(A)
$$
defined as follows. For two elements $a'_0[a'_1|a'_2|\ldots
|a'_n],a''_0[a''_1|a''_2|\ldots |a''_m]\in \rC_\bullet(A)$ the
shuffle product is given by the formula:
\begin{equation}\label{shuff}
\sh(a'_0[a'_1|a'_2|\ldots |a'_n]\otimes
a''_0[a''_1|a''_2|\ldots
|a''_m])=(-1)^{\ast}\cdot a'_0a''_0\,\sh_{nm}[a'_1|\ldots|a'_n| a''_1|\ldots|a''_m]
\end{equation}
Here $\ast=|a''_0|(|sa'_1|+\ldots+|sa'_n|)$ and
$$
\sh_{nm}[x_1|\ldots|x_n|x_{n+1}|\ldots|x_{n+m}]=\sum_{\sigma}\pm[x_{\sigma^{-1}(1)}|\ldots|
x_{\sigma^{-1}(n)}|
x_{\sigma^{-1}(n+1)}|\ldots|x_{\sigma^{-1}(n+m)}]
$$
where the sum is taken over all permutations that don't shuffle the
first $n$ and the last $m$ elements and the sign in front of each
summand is computed according to the following rule: for two
homogeneous elements $x,y$, the transposition
$[\,\ldots|x|y|\ldots\,]\to[\,\ldots|y|x|\ldots\,]$ contributes
$(-1)^{|x||y|}$ to the sign.

Now let $B$ be another DG algebra. Denote by $\iota^A,\iota^B$ the
natural embeddings $$A\to A\otimes B, \quad B\to A\otimes B.$$ They
induce morphisms of complexes:
$$
\rC(\iota^A): \rC_\bullet(A)\to \rC_\bullet(A\otimes B),\quad
\rC(\iota^B): \rC_\bullet(B)\to \rC_\bullet(A\otimes B).
$$

\medskip
\begin{theorem}\label{kunneth}
The composition $\K$ of the maps
\[
  \begin{CD}
\rC_\bullet(A)\otimes \rC_\bullet(B)
@>\rC(\iota^A)\otimes\rC(\iota^B)>>  \rC_\bullet(A\otimes B)\otimes
\rC_\bullet(A\otimes B) @>\sh>>  \rC_\bullet(A\otimes B)
  \end{CD}
\]
respects the differentials and induces a quasi-isomorphism of
complexes.
\end{theorem}
\medskip
The morphism $\K:\rC_\bullet(A)\otimes \rC_\bullet(B)\to
\rC_\bullet(A\otimes B)$ defined above admits a generalization to
the case of DG categories. Namely, let $\cA$ and $\cB$ be two
(small) DG categories. Fix a set $\{X_0,X_1,\ldots,X_n\}$ of objects
of $\cA$ and a set $\{Y_0,Y_1,\ldots,Y_m\}$ of objects of $\cB$. For
two elements
$$
f_n[f_{n-1}|\ldots |f_0]\in\Hom_{\cA}(X_n,X_0)\otimes
\Hom_{\cA}(X_{n-1},X_n)[1]\otimes\ldots\otimes
\Hom_{\cA}(X_0,X_1)[1],
$$
$$
g_m[g_{m-1}|\ldots |g_0]\in\Hom_{\cB}(Y_m,Y_0)\otimes
\Hom_{\cB}(Y_{m-1},Y_m)[1]\otimes\ldots\otimes
\Hom_{\cB}(Y_0,Y_1)[1]
$$
define $\K\left(f_n[f_{n-1}|f_{n-2}|\ldots |f_0]\otimes
g_m[g_{m-1}|g_{m-2}|\ldots |g_0]\right)$ as
$$
\pm(f_n\otimes g_m)\,\sh_{nm}[f_{n-1}|\ldots |f_0|g_{m-1}|\ldots|g_0],
$$
where the sign is computed as before and $\sh_{nm}$ is defined by
the formula
$$
[f_{n-1}\otimes 1_{Y_m}|\ldots |f_0\otimes
1_{Y_m}|1_{X_0}\otimes g_{m-1}|\ldots|1_{X_0}\otimes g_0]+
$$
$$
+(-1)^{|sf_0||sg_{m-1}|}[f_{n-1}\otimes 1_{Y_m}|\ldots
|1_{X_1}\otimes g_{m-1}|f_0\otimes
1_{Y_{m-1}}|\ldots|1_{X_0}\otimes g_0]+\ldots
$$
Other terms in this sum are obtained from the first one by shuffling
the $f$-terms with the $g$-terms according to the following rule:
$$
[\,\ldots|f_{k}\otimes 1_{Y_{l+1}}|1_{X_k}\otimes
g_{l}|\ldots\,]\to(-1)^{|sf_k||sg_l|}[\,\ldots|1_{X_{k+1}}\otimes
g_{l}|f_{k}\otimes 1_{Y_{l}}|\ldots\,]
$$

\medskip

Let $A$ and $B$ be two DG algebras. We have the obvious embedding of
DG categories
$$
\per A \otimes \per B \to \per(A\otimes B),
$$
which induces a morphism of complexes
$$
\rC_\bullet(\per A \otimes \per B)\to\rC_\bullet(\per (A\otimes B)).
$$
Let us denote the composition
\begin{equation}\label{nk}
\rC_\bullet(\per A)\otimes \rC_\bullet(\per B)\stackrel{\K}\to\rC_\bullet(\per A \otimes \per B)
\to\rC_\bullet(\per (A\otimes B))
\end{equation}
by the same letter $\K$. As an immediate corollary of Theorems
\ref{hhper} and \ref{kunneth}, we get the following result:

\medskip
\begin{proposition}
The map $\K: \rC_\bullet(\per A)\otimes \rC_\bullet(\per
B)\to\rC_\bullet(\per (A\otimes B))$ is a quasi-isomorphism.
\end{proposition}
\noindent Indeed, we have the commutative diagram
\begin{displaymath}
\xymatrix{ \rC_\bullet(\per A)\otimes \rC_\bullet(\per B) \ar[r] &
\rC_\bullet(\per(A\otimes B))\\
\rC_\bullet(A)\otimes \rC_\bullet(B) \ar[u]\ar[r] &
\rC_\bullet(A\otimes B)\ar[u] }
\end{displaymath}
in which the vertical arrows and the arrow on the bottom are
quasi-isomorphisms.

\medskip

Finally, we will formulate two more results about the K\"{u}nneth
map (\ref{nk}). Both results follow directly from the definition of
$\K$.

\medskip
\begin{proposition}\label{tens}
Let $A$, $B$, and $C$ be three DG algebras. The diagram
\begin{displaymath}
\xymatrix{ \rC_\bullet(\per A)\otimes \rC_\bullet(\per B) \otimes
\rC_\bullet(\per C) \ar[r]^{\,\,\,\K\otimes 1}\ar[d]_{1\otimes \K} &
\rC_\bullet(\per (A\otimes B)) \otimes
\rC_\bullet(\per C) \ar[d]^{\K}\\
\rC_\bullet(\per A)\otimes \rC_\bullet(\per (B\otimes C))
\ar[r]^{\K} & \rC_\bullet(\per( A\otimes B\otimes C))}
\end{displaymath}
commutes. In other words, $\K$ is associative.
\end{proposition}

\medskip

\begin{proposition}\label{tensfun}
Let $A,B,C,D$ be DG algebras. Let $X\in\Mod(A^\op\otimes C)$ and
$Y\in\Mod(B^\op\otimes D)$ be bimodules satisfying the conditions of
Corollary \ref{cor1}. Then the diagram
\begin{displaymath}
\xymatrix{ \rC_\bullet(\per A)\otimes \rC_\bullet(\per B)
\ar[r]^{\,\,\,\K}\ar[d]_{\rC(T_X)\otimes \rC(T_Y)} &
\rC_\bullet(\per (A\otimes B)) \ar[d]^{\rC(T_{X\otimes_k Y})}\\
\rC_\bullet(\per C)\otimes \rC_\bullet(\per D) \ar[r]^{\,\,\,\K} &
\rC_\bullet(\per (C\otimes D))}
\end{displaymath}
commutes.
\end{proposition}

\medskip
\section{Hirzebruch-Riemann-Roch theorem}
\subsection{Euler character}\label{ec} Let $A$ be a DG algebra and $N$ a
perfect right DG $A$-module. Consider the DG functor $T_N=-\otimes_k
N: \per k\to\per A$. The Euler class $\eu(N)\in\rHH_0(\per A)$ is
defined by the formula (cf. \cite{BNT},\cite{K1})
$$
\eu(N)=\rHH(T_N)(1).
$$
In other words, $\eu(N)$ is the class of the identity morphism $1_N$
in $\rHH_0(\per A)$.

Let us list some basic properties of the Euler character map.

The following statement follows from Theorem \ref{hominv}:
\medskip
\begin{proposition}\label{p1}
If $N,M\in\per A$ are homotopically equivalent then $\eu(N)=\eu(M)$.
In other words, $\eu$ descends to objects of $\Ho(\per A)$.
\end{proposition}
\medskip
The following result means that the Euler class descends to the
Grothendieck group of the triangulated category $\Ho(\per A)$.
\medskip
\begin{proposition}\label{p2}
For any $N\in\per A$ one has $\eu(N[1])=-\eu(N)$ and for any
triangle $L\stackrel{p}\to M\stackrel{q}\to N \stackrel{r}\to L[1]$
in $\Ho(\per A)$ one has
\begin{equation}\label{add}
\eu(M)=\eu(L)+\eu(N).
\end{equation}
\end{proposition}
Let us prove the first part. We have to show that $1_N+1_{N[1]}$ is
homologous to 0 in $\rC_\bullet(\per A)$. Denote by $1_{N,N[1]}$
(resp. $1_{N[1],N}$) the identity endomorphism of $N$ viewed as a
morphism from $N$ to $N[1]$ (resp. from $N[1]$ to $N$). Then
\begin{eqnarray*}
b(1_{N,N[1]}[1_{N[1],N}])=b_1(1_{N,N[1]}[1_{N[1],N}])=\\
=-(1_{N,N[1]}1_{N[1],N}+1_{N[1],N}1_{N,N[1]})=-(1_{N[1]}+1_N)
\end{eqnarray*}

Let us prove the second part. By Proposition \ref{p1}, it suffices
to prove (\ref{add}) for $N=\cone(p)$. Consider the following
morphisms:
$$
j_1=\begin{pmatrix} 1_{L[1]} \\ 0 \end{pmatrix}:
L[1]\to\cone(p),\quad q_1=\begin{pmatrix}1_{L[1]} & 0 \end{pmatrix}:
\cone(p)\to L[1],
$$
$$
j_2=\begin{pmatrix}0 \\1_{M} \end{pmatrix}: M\to\cone(p),\quad
q_2=\begin{pmatrix}0 & 1_{M} \end{pmatrix}: \cone(p)\to M.
$$
It is easy to see that
$$
d(j_1)=j_2\cdot p,\quad d(q_1)=0,\quad d(j_2)=0,\quad d(q_2)=-p\cdot q_1.
$$
(In these formulas, $p$ is viewed as a degree 1 morphism from $L[1]$
to $M$.) The following computation finishes the proof:
\begin{eqnarray*}
1_{\cone(p)}-1_{L[1]}-1_{M}=j_1q_1+j_2q_2-q_1j_1-q_2j_2=[j_1,q_1]+[j_2,q_2]\\
=b(j_1[q_1]+j_2[q_2])-b_0(j_1[q_1]+j_2[q_2])=b(j_1[q_1]+j_2[q_2])-(d(j_1)[q_1]-j_2[d(q_2)])\\
=b(j_1[q_1]+j_2[q_2])-(j_2p[q_1]+j_2[pq_1])=b(j_1[q_1]+j_2[q_2]-j_2[p|q_1]).
\end{eqnarray*}
\medskip

To formulate the main result of this section, we need a pairing
$$\rHH_n(\per A)\times\rHH_{-n}(\per A^\op)\to k, \quad
n\in\Z,$$ where $A$ is a proper DG algebra. Here is the definition.

Let us equip $A$ with a left DG $A\otimes A^\op$-module structure as
follows:
$$
(a'\otimes a'')a=(-1)^{|a''||a|}a'aa''.
$$
We will denote the resulting $A$-bimodule by $\Delta$.

Consider the DG functor:
$$
T_\Delta: \Mod(A\otimes A^\op)\to\Mod k, \quad N\mapsto
N\otimes_{A\otimes A^\op} A
$$

The following proposition is an immediately consequence of Corollary
\ref{cor1}.

\medskip
\begin{proposition}
If $A$ is proper then $T_\Delta$ induces a DG functor $\per(A\otimes
A^\op)\to\per k$.
\end{proposition}
\medskip

We can use this to define a pairing
\begin{equation}\label{paironperf}
\langle\,,\,\rangle: \rHH_n(\per A)\times\rHH_{-n}(\per A^\op)\to
k, \quad n\in\Z
\end{equation}
via the composition of morphisms of complexes
\[
  \begin{CD}
  \rC_\bullet(\per A)\otimes
\rC_\bullet(\per A^\op)@>\K>> \rC_\bullet(\per(A\otimes
A^\op))@>\rC(T_\Delta)>> C_\bullet(\per k)
  \end{CD}
\]
and the fact that $\rHH_n(\per k)\simeq\rHH_n(k)$ is $k$, if $n=0$,
and $0$ otherwise.

Before we formulate the main result of this section, let us
introduce the following notation. For a bimodule
$X\in\per(A^\op\otimes B)$ we will denote by $\eu'(X)$ the element
$$\K^{-1}(\eu(X))\in\bigoplus_n\rHH_{-n}(\per A^\op)\otimes
\rHH_{n}(\per B),$$ where $\K$ is the K\"{u}nneth isomorphism.

\medskip
\begin{theorem}\label{main}
Let $A$ be a proper DG algebra, $B$ an arbitrary DG algebra, and $X$
any object of $\per(A^\op\otimes B)$. If $y\in\rHH_\bullet(\per A)$
then $ \rHH(T_X)(y)=\langle\,y\,,\,\eu'(X)\,\rangle. $ That is, if
$$
\eu'(X)=\sum_nx'_{-n}\otimes x''_{n}\in\bigoplus_n\rHH_{-n}(\per
A^\op)\otimes\rHH_{n}(\per B),$$ then $
\rHH(T_X)(y)=\sum_n\langle\,y\,,\,x'_{-n}\,\rangle\cdot x''_n. $
\end{theorem}

To prove this, observe that $T_X$ can be described as a composition
of the following DG functors
\[
  \begin{CD}
\per A@>-\otimes_k X>>\per(A\otimes A^\op\otimes
B)@>T_{\Delta\otimes_k B}>>\per B
  \end{CD}
\]
Thus, $\rHH(T_X)=\rHH(T_{\Delta\otimes_k B})\circ\rHH(-\otimes_k
X)$. It follows from the definition of the K\"{u}nneth isomorphism
$\K$ that the diagram
\begin{displaymath}
\xymatrix{ \rHH_\bullet(\per
A)\ar[d]_{1\otimes\eu(X)}\ar[r]^{\rHH(-\otimes_k X)} &
\rHH_\bullet(\per(A\otimes A^\op\otimes B))\\
\rHH_\bullet(\per A)\otimes \rHH_0(\per(A^\op\otimes B))\ar[ur]^{\K}
&}
\end{displaymath}
commutes. By conjugating with $1\otimes\K$, we get the following
commutative diagram:
\begin{displaymath}
\xymatrix{ \rHH_\bullet(\per
A)\ar[d]_{1\otimes\eu'(X)}\ar[r]^{\rHH(-\otimes_k X)} &
\rHH_\bullet(\per(A\otimes A^\op\otimes B))\\
\rHH_\bullet(\per A)\otimes \rHH_\bullet(\per A^\op)\otimes
\rHH_\bullet(\per B)\ar[ur]^{\K\circ(1\otimes\K)} &}
\end{displaymath}

Furthermore, by Proposition \ref{tensfun} the diagram
\begin{displaymath}
\xymatrix{\rHH_\bullet(\per(A\otimes A^\op\otimes
B))\ar[d]_{\K^{-1}}\ar[r]^{\rHH(T_{\Delta\otimes_k B})} &
\,\,\,\,\,\,\rHH_\bullet(\per(k\otimes B))\simeq\rHH_\bullet(\per B)\\
\rHH_\bullet(\per(A\otimes A^\op))\otimes \rHH_\bullet(\per
B)\ar[r]^{\quad\rHH(T_\Delta)\otimes1} & \rHH_\bullet(\per k)\otimes
\rHH_\bullet(\per B)\ar[u]^\K}
\end{displaymath}
commutes. Conjugating with $\K\otimes1$ gives us the following
commutative diagram:
\begin{displaymath}
\xymatrix{\rHH_\bullet(\per(A\otimes A^\op\otimes
B))\ar[d]_{(\K^{-1}\otimes1)\K^{-1}}\ar[r]^{\rHH(T_{\Delta\otimes_k
B})} &
\rHH_\bullet(\per(k\otimes B))\simeq\rHH_\bullet(\per B)\\
\rHH_\bullet(\per A)\otimes\rHH_\bullet(\per A^\op)\otimes
\rHH_\bullet(\per B)\ar[r]^{\qquad\qquad(\rHH(T_\Delta)\K)\otimes1}
& \,\,\,\,\,\rHH_\bullet(\per k)\otimes \rHH_\bullet(\per
B)\ar[u]^\K}
\end{displaymath}
By concatenating the top arrows of the former and the latter
diagrams, we get the following result:
$$
\rHH(T_{\Delta\otimes_k B})\circ\rHH(-\otimes_k
X)=\K\circ((\rHH(T_\Delta)
\K)\otimes1)\circ(\K^{-1}\otimes1)\circ\K^{-1}\circ\K\circ(1\otimes\K)\circ(1\otimes\eu'(X)).
$$
By associativity of the K\"{u}nneth isomorphism (Proposition
\ref{tens}), the latter product is nothing but
$\K\circ((\rHH(T_\Delta)\K)\otimes1)\circ(1\otimes\eu'(X))$ which
finishes the proof.

\medskip

Theorem \ref{main} generates several corollaries. The first one, the
Hirzebruch-Riemann-Roch type formula, will be formulated and proved
in the next section. Another corollary, which concerns the so-called
smooth DG algebras, will be described in Section \ref{smooth}.

\medskip

\subsection{Hirzebruch-Riemann-Roch theorem}\label{hrr}

Essentially, the Hirzebruch-Riemann-Roch theorem is the following
result:

\medskip

\begin{theorem}\label{rr1}
Let $A$ be a proper DG algebra. Then, for any $N\in\per A$,
$M\in\per A^\op$,
\begin{equation}\label{rrh}
\sum_n(-1)^n\dim \rH^n(N\otimes_A M)=\langle \eu(N),\eu(M)\rangle.
\end{equation}
\end{theorem}
This theorem is an easy corollary of the results of the previous
section. Indeed, consider the DG functors:
$$
T_N=-\otimes_k N:\per k\to\per A,\quad T_M=-\otimes_A M:\per A\to\per
k,$$$$T_{N\otimes_AM} =-\otimes_k (N\otimes_AM):\per k\to\per k.
$$
Clearly, $T_{N\otimes_AM}=T_{M}T_{N}$ and, by Theorem \ref{main}, we
get the equality
$$
\eu(N\otimes_AM)=\langle\eu(N),\eu(M)\rangle.
$$
What remains is to observe that, for a perfect DG $k$-module $X$,
$$
\eu(X)=\sum_n(-1)^n\dim \rH^n(X).
$$
This latter statement is a corollary of Propositions \ref{p1} and
\ref{p2}, along with the fact that $X$ is homotopy equivalent to
$\rH^\bullet(X)$.

\medskip

Let us explain how one can compute the right-hand side of
(\ref{rrh}).

First of all, observe that, by Theorem \ref{hhper}, the pairing
(\ref{paironperf}) induces a pairing on
$\rHH_\bullet(A)\times\rHH_\bullet(A^\op)$. Let us fix two cycles
$$\sum_{a} a_0[a_1|\ldots |a_l]\in\rC_\bullet(A),\quad \sum_{b}
b_0[b_1|\ldots |b_m]\in\rC_\bullet(A^\op)$$ ($\sum$ indicates that
$a$ and $b$ are sums of several terms) and denote by $a$ (resp. $b$)
the corresponding elements in $\rHH_\bullet(A)$ (resp.
$\rHH_\bullet(A^\op)$). Let us describe $\langle \,a,b\,\rangle$
more explicitly.

Consider the composition of DG functors $$A\otimes A^\op\to
\per(A\otimes A^\op)\stackrel{T_\Delta}\to\per k,$$ where $A\otimes
A^\op$ is viewed as a DG category with one object. Clearly, the
unique object of $A\otimes A^\op$ gets mapped under this composition
to $A\in\per k$ and an element $x\otimes y\in A\otimes A^\op$,
viewed as a morphism in the DG category $A\otimes A^\op$, gets
mapped to the operator $L(x)R(y)\in\End_k(A)$, where
$$
L(x): c\mapsto xc,\quad R(y): c\mapsto (-1)^{|c||y|}cy
$$
are the operators of left multiplication with $x$ resp. right
multiplication with $y$.

Since the operators of left multiplication commute with operators of
right multiplication, we can define a product
\begin{eqnarray}\label{wedge}
a\wedge b=\sum_{a,b}\pm
L(a_0)R(b_0)\sh_{lm}[L(a_1)|\ldots |L(a_l)|R(b_1)|\ldots
|R(b_m)]
\end{eqnarray}
on $\rHH_\bullet(A)\times\rHH_\bullet(A^\op)$ with values in
$\rHH_\bullet(\End_k(A))$ (the formula for $\pm$ and the definition
of $\sh_{lm}$ are the same as in (\ref{shuff})). Then
\begin{eqnarray}\label{pairing}
\langle \,a,b\,\rangle=\int a\wedge b
\end{eqnarray}
where $\int$ is defined as follows. Let $X$ be a perfect DG
$k$-module. Then we have an embedding of DG categories\footnote{For
a complex $X$ of vector spaces $\End_k(X)$ stands for the DG algebra
$\oplus_n\End^n_k(X)$ where $\End^n_k(X)$ is the subspace of degree
$n$ linear maps.} $\End_k(X)\to\per k$ which sends the unique object
of the first category to $X$, viewed as an object of $\per k$. Then
$\int$ is the map from $\rHH_\bullet(\End_k(X))$ to
$\rHH_\bullet(\per k)\simeq k$ induced by this embedding.

Furthermore, let us use the notation $\eul(N)$ to denote the element
in $\rHH_0(A)$ corresponding to $\eu(N)$ under the isomorphism
$\rHH_0(A)\to\rHH_0(\per A)$. We are ready to rewrite the right-hand
side of (\ref{rrh}):
$$
\langle \eu(N),\eu(M)\rangle=\int\eul(N)\wedge\eul(M).
$$
It turns out that there are very explicit formulas for $\int$ and
$\eul$ which will be derived in the next section.

\medskip

To conclude this section, we will rewrite the
Hirzebruch-Riemann-Roch formula in a more conventional form. Namely,
we will use (\ref{rrh}) to derive a formula that expresses the Euler
form
$$
\chi(M,N)=\sum_n(-1)^n\dim\Hom_{\Ho(\per A)}(M,N[n])
$$
in terms of the Euler classes of $M$ and $N$, where $M$ and $N$ are
two perfect DG $A$-modules.

Consider the following (contravariant) DG functor
\begin{equation}\label{qe}
D: \Mod A\to\Mod A^\op, \quad M\mapsto DM=\Hom_{\Mod A}(M,
A).
\end{equation}
It is not hard to show that this DG functor preserves perfect
modules. Moreover, its square is isomorphic to the identity
endofunctor of $\per A$ and, thus, $D$ is a quasi-equivalence of the
DG categories $(\per A)^\op$ and $\per A^\op$. The crucial property
of this functor is the following fact: for any perfect DG
$A$-modules there is a natural quasi-isomorphism of complexes
$$
N\otimes_A DM\cong\Hom_{\per A}(M,N).
$$

Thus, the formula (\ref{rrh}) can be written as follows: for any
$N,M\in\per A$
\begin{equation}\label{rrh1}
\chi(M,N)=\langle \eu(N),\eu(DM)\rangle=\int\eul(N)\wedge\eul(DM).
\end{equation}

Finally, we notice that $\eu(DM)$ (and $\eul(DM)$) can be expressed
in terms of $\eu(M)$ (resp. $\eul(M)$). This is based on the
following result (see Appendix \ref{proof}):

\medskip
\begin{proposition}\label{veeiso}
For any DG algebra, the formula
\begin{equation}\label{aaop}
(a_0[a_1|a_2|\ldots |a_n])^\vee=(-1)^{n+\sum_{1\leq i<j\leq
n}|sa_i||sa_j|}a_0[a_n|a_{n-1}|\ldots |a_1].
\end{equation}
defines a quasi-isomorphism $\,^\vee: \rC_\bullet(A)\to
\rC_\bullet(A^\op)$.
\end{proposition}
\medskip

One can generalize the above formulas to the case of an arbitrary DG
category to get a quasi-isomorphism $\,^\vee: \rC_\bullet(\cA)\to
\rC_\bullet(\cA^\op)$. In the case $\cA=\per A$ one can compose it
with $\rC(D): \rC_\bullet((\per A)^\op)\to \rC_\bullet(\per A^\op)$
to get a quasi-isomorphism $\,^\vee: \rC_\bullet(\per A)\to
\rC_\bullet(\per A^\op)$. It is immediate that
$\eu(DM)=\eu(M)^\vee$. It is also true, but is less obvious, that
$\eul(DM)=\eul(M)^\vee$. This latter observation follows from the
fact that the two quasi-isomorphisms $\,^\vee: \rC_\bullet(A)\to
\rC_\bullet(A^\op)$ and $\,^\vee: \rC_\bullet(\cA)\to
\rC_\bullet(\cA^\op)$ agree under the embeddings $A\to\per A$ and
$A^\op\to\per A^\op$.

\medskip

So here is the Hirzebruch-Riemann-Roch formula in its ultimate form:
\begin{equation}\label{rrh2}
\chi(M,N)=\langle \eu(N),\eu(M)^\vee\rangle=\int\eul(N)\wedge\eul(M)^\vee.
\end{equation}

\medskip
\section{On the computational aspect of the HRR theorem}
\subsection{Computing Euler classes}\label{cec}
The aim of this section is to explain how to compute the Euler class
$\eul(N)\in\rHH_0(A)$ of a perfect DG $A$-module.

The definition of a finitely generated semi-free module we gave in
Section \ref{perf} is convenient for proving theorems but it is not
explicit enough for the purposes of this section. A more explicit
description was given in \cite{BK} and we will begin by recalling
it.

Let $A$ be a DG algebra. Let $\free A$ be the DG subcategory in
$\per A$ whose objects are finitely generated free DG $A$-modules,
i.e. direct sums of modules of the form $$A[r]=k[r]\otimes A, \quad
r\in\Z.$$ Clearly, $$\Hom_{\free A}(A[r], A[s])=\Hom_{\per A}(A[r],
A[s])\simeq A[s-r].$$ The differential on the morphism spaces of the
DG category $\free A$, as well as on the free modules themselves,
will be denoted by $d_{\free}$.

The alternative description of finitely generated semi-free modules
is based on the notion of a twisted $A$-modules. These are objects
of a larger DG subcategory $\Tw A\supset\free A$ in $\per A$.
Namely, a twisted $A$-module is a right DG $A$-module of the form
$(\bigoplus\limits_{j=1}^nA[r_j], d_{\free}+\alpha)$, where
$\alpha=(\alpha_{ij})$ is a strictly upper triangular $n\times
n$-matrix of morphisms $\alpha_{ij}\in \Hom_{\free A}^1(A[r_j],
A[r_i])$ satisfying the Maurer-Cartan equation $$
d_{\free}(\alpha)+\alpha\cdot\alpha=0.$$ Clearly, the differential
$d_{\Tw}$ on $\Hom_{\per A}((\bigoplus\limits_{j=1}^nA[r_j],
d_{\free }+\alpha), (\bigoplus\limits_{i=1}^mA[s_i],
d_{\free}+\beta))$ is given by the formula $$
d_{\Tw}(f)=d_{\free}(f)+\beta\cdot f-(-1)^{|f|}f\cdot\alpha. $$ It
is not hard to show that any finitely generated semi-free module is
isomorphic to a twisted $A$-modules.

The main result of this section is a formula for the Euler class of
a homotopy direct summand of a twisted $A$-module. Its formulation
involves a (super-)trace map\footnote{In the case of an associative
algebra, this map coincides with the well-known trace map from
Section 1.2 of \cite{L}.} $\str$ which we will describe now.

Let $N$ be a DG $A$-module which is isomorphic to
$\bigoplus\limits_{j=1}^nA[r_j]$ as a graded $A$-module. Fix $m$
homogeneous endomorphisms of $N$:
$$A',A'',\ldots, A^{(m)}\in\End_{\per A}(N).$$
Thus, each $A^{(k)}$ is an $n\times n$-matrix $(e(r_i,r_j)\otimes
a^{(k)}_{ij})$ of morphisms $$e(r_i,r_j)\otimes a^{(k)}_{ij} \in
\Hom_{\per A}(A[r_j], A[r_i]),$$ where $a^{(k)}_{ij}\in A$ and
$e(r_i,r_j)\in\Hom_{\per A}(A[r_j], A[r_i])$ is the morphism that
sends the generator of $A[r_j]$ to the generator of $A[r_i]$. The
endomorphisms give rise to an element $A'[A''|\ldots|A^{(m)}]$ of
the Hochschild chain complex of the DG category $\per A$. Let us
define $\str(A'[A''|\ldots|A^{(m)}])\in\rC_\bullet(A)$ by the
formula
\begin{eqnarray*}
\str(A'[A''|\ldots|A^{(m)}])=\sum_{j=1}^n\sum_{i_1,i_2,\ldots,i_{m-1}}
(-1)^{\ast}
\cdot a'_{ji_1}[a''_{i_1i_2}|\ldots|a^{(m)}_{i_{m-1}j}],
\end{eqnarray*}
where $ \ast=
r_{i_1}+(r_{i_1}-r_{j})|a'_{ji_1}|+(r_{i_2}-r_{j})|sa''_{i_1i_2}|+\ldots+(r_{i_{m-1}}-r_{j})|
sa^{(m-1)}_{i_{m-2}i_{m-1}}|.$

\medskip
\begin{theorem}\label{eul}
Let $N_{\alpha}=(\bigoplus\limits_{j=1}^nA[r_j], d_{\free}+\alpha)$
and $L$ be a homotopy direct summand of $N_{\alpha}$ corresponding
to a homotopy idempotent $\pi:N_{\alpha}\to N_{\alpha}$. Then
$$
\eul(L)=\sum_{l=0}^{n-1}
(-1)^l\str(\pi[\underbrace{\alpha|\ldots|\alpha}_{l}])
$$
\end{theorem}

\medskip

Let us prove the theorem.
\medskip
\begin{lemma}
In the above notation, $\eu(L)=\pi$.
\end{lemma}
We have to show that $1_{L}\in\End^0_{\per A}(L)$ and
$\pi\in\End^0_{\per A}(N_{\alpha})$ define the same element of
$\rHH_0(\per A)$. Let us fix some degree 0 closed morphisms
$f:N_{\alpha}\to L$ and $g:L\to N_{\alpha}$ such that
$$fg=1_{L}+[d_L,H_L],\quad gf=\pi+[d_{N_{\alpha}},H_{N_{\alpha}}]$$ (see
Section \ref{1}). Then
$$
1_{L}-\pi=b(f[g]+H_{N_{\alpha}}-H_L).
$$
The lemma is proved.

\medskip

Let $N_{\alpha}=(\bigoplus\limits_{j=1}^nA[r_j], d_{\free}+\alpha)$
and $\pi$ be as before. Let us introduce some new notations. We will
write $N_0$ to denote the free DG $A$-module
$(\bigoplus\limits_{j=1}^nA[r_j], d_{\free})$. For an endomorphism
$f\in\End_{\per A}(N_{\alpha})$, $\widetilde{f}$ (resp.
$\overrightarrow{f}$, $\overleftarrow{f}$) will stand for $f$ viewed
as an element of $\End_{\per A}(N_0)$ (resp. $\Hom_{\per
A}(N_{\alpha},N_0)$, $\Hom_{\per A}(N_0,N_{\alpha})$). For a
morphism $g$ we will write $g_{ij}$ (resp. $g_{i\ast}$, $g_{\ast
j}$) for the $n\times n$-matrix, viewed as a morphism between the
same modules, whose $ij$-th entry (resp. $i$-th row, $j$-th column)
coincides with that of $g$ and other entries (resp. rows, columns)
are 0.

The following lemma is a straightforward consequence of the
definition of $\str$:

\medskip
\begin{lemma}\label{l0}
One has
$$
\sum_{l=0}^{n-1}(-1)^l
\str({\pi[\underbrace{\alpha|\ldots|\alpha}_l]})=\sum_{l=0}^{n-1}\sum_{i_0,i_1,\ldots,i_l}
(-1)^l \str({\widetilde{\pi}_{i_0i_1}[\widetilde{\alpha}_{i_1i_2}|\ldots|\widetilde{\alpha}_{i_li_0}]})
$$
in $\rHH_0(\per A)$.
\end{lemma}

\medskip
The next lemma is less straightforward:
\medskip

\begin{lemma}\label{l1}
One has
$$
\pi=\sum_{l=0}^{n-1}\sum_{i_0,i_1,\ldots,i_l}
(-1)^l {\widetilde{\pi}_{i_0i_1}[\widetilde{\alpha}_{i_1i_2}|\ldots|\widetilde{\alpha}_{i_li_0}]}
$$
in $\rHH_0(\per A)$.
\end{lemma}
To prove this, pick a large $N$ and apply the differential $b$ to
the element
$$
\sum_{l=0}^{N}\sum_{i_0,i_1,\ldots,i_l}
(-1)^l {\overrightarrow{\pi}_{i_0\ast}[\overleftarrow{1}_{\ast i_1}|\widetilde{\alpha}_{i_1i_2}|\ldots|\widetilde{\alpha}_{i_li_0}]}.
$$
Let us begin by computing the $b_0$-component:
\begin{eqnarray*}
b_0({\overrightarrow{\pi}_{i_0\ast}[\overleftarrow{1}_{\ast i_1}|\widetilde{\alpha}_{i_1i_2}|\ldots|\widetilde{\alpha}_{i_li_0}]})=
d_{\Tw}({\overrightarrow{\pi}_{i_0\ast})[\overleftarrow{1}_{\ast i_1}|\widetilde{\alpha}_{i_1i_2}|\ldots|\widetilde{\alpha}_{i_li_0}]}\\
-{\overrightarrow{\pi}_{i_0\ast}[d_{\Tw}(\overleftarrow{1}_{\ast i_1})|\widetilde{\alpha}_{i_1i_2}|\ldots|\widetilde{\alpha}_{i_li_0}]}\\
+\sum_{m=1}^{l}
{\overrightarrow{\pi}_{i_0\ast}[\overleftarrow{1}_{\ast i_1}|\widetilde{\alpha}_{i_1i_2}|\ldots|d_{\free}(\widetilde{\alpha}_{i_mi_{m+1}})
|\ldots|\widetilde{\alpha}_{i_li_0}]}.
\end{eqnarray*}
Recall that $\pi$ is closed, i.e.
$d_{\free}(\pi)+\alpha\pi-\pi\alpha=0$. Therefore
\begin{eqnarray*}
d_{\Tw}(\overrightarrow{\pi}_{i_0\ast})=d_{\free}(\overrightarrow{\pi}_{i_0\ast})
-\overrightarrow{\pi}_{i_0\ast}\alpha
=(\overrightarrow{\pi}\alpha)_{i_0\ast}-(\widetilde{\alpha}
\overrightarrow{\pi})_{i_0\ast}-\overrightarrow{\pi}_{i_0\ast}\alpha\\=
-(\widetilde{\alpha}\overrightarrow{\pi})_{i_0\ast}=-\sum_{k=1}^n
\widetilde{\alpha}_{i_0k}\overrightarrow{\pi}_{k\ast}.
\end{eqnarray*}
Furthermore,
$$d_{\Tw}(\overleftarrow{1}_{\ast i_1})=\alpha\overleftarrow{1}_{\ast i_1}=
\overleftarrow{\alpha}_{\ast i_1},\quad  d_{\free}(\widetilde{\alpha}_{i_mi_{m+1}})=-\sum_{k=1}^n\widetilde{\alpha}_{i_mk}
\widetilde{\alpha}_{ki_{m+1}}.$$
Thus,
\begin{eqnarray*}\label{b0comp}
b_0({\overrightarrow{\pi}_{i_0\ast}[\overleftarrow{1}_{\ast i_1}|\widetilde{\alpha}_{i_1i_2}|\ldots|\widetilde{\alpha}_{i_li_0}]})=-\sum_{k=1}^n
\widetilde{\alpha}_{i_0k}{\overrightarrow{\pi}_{k\ast}[\overleftarrow{1}_{\ast i_1}|\widetilde{\alpha}_{i_1i_2}|\ldots|\widetilde{\alpha}_{i_li_0}]}\\
-\overrightarrow{\pi}_{i_0\ast}[\overleftarrow{\alpha}_{\ast i_1}|\widetilde{\alpha}_{i_1i_2}|\ldots|\widetilde{\alpha}_{i_li_0}]\\
-\sum_{m=1}^{l}\sum_{k=1}^n
{\overrightarrow{\pi}_{i_0\ast}[\overleftarrow{1}_{\ast i_1}|\widetilde{\alpha}_{i_1i_2}|\ldots|\widetilde{\alpha}_{i_mk}
\widetilde{\alpha}_{ki_{m+1}}|\ldots|\widetilde{\alpha}_{i_li_0}]}.
\end{eqnarray*}

Let us compute now the $b_1$-component. Clearly,
$b_1(\overrightarrow{\pi}_{i_0\ast} [\overleftarrow{1}_{\ast
i_0}])=\widetilde{\pi}_{i_0i_0}-\pi_{i_0i_0}$ and for $l\geq1$
\begin{eqnarray*}
b_1({\overrightarrow{\pi}_{i_0\ast}[\overleftarrow{1}_{\ast i_1}|\widetilde{\alpha}_{i_1i_2}|\ldots|\widetilde{\alpha}_{i_li_0}]})
={\widetilde{\pi}_{i_0i_1}[\widetilde{\alpha}_{i_1i_2}|\ldots|\widetilde{\alpha}_{i_li_0}]}
-{\overrightarrow{\pi}_{i_0\ast}[\overleftarrow{\alpha}_{i_1i_2}|\ldots|
\widetilde{\alpha}_{i_li_0}]}
\\
-\sum_{m=1}^{l-1}{\overrightarrow{\pi}_{i_0\ast}[\overleftarrow{1}_{\ast i_1}|\widetilde{\alpha}_{i_1i_2}|\ldots|\widetilde{\alpha}_{i_mi_{m+1}}
\widetilde{\alpha}_{i_{m+1}i_{m+2}}|\ldots|\widetilde{\alpha}_{i_li_0}]}\\
-{\widetilde{\alpha}_{i_li_0}\overrightarrow{\pi}_{i_0\ast}[\overleftarrow{1}_{\ast i_1}|\widetilde{\alpha}_{i_1i_2}|\ldots|\widetilde{\alpha}_{i_{l-1}i_l}]}.
\end{eqnarray*}
To finish the proof of Lemma \ref{l1}, one needs to add the results
of the above two computations, take the sum over $l$ and $i_0,
i_1,\ldots,i_l$, and observe that the right-hand side of the formula
for $b_0({\overrightarrow{\pi}_{i_0\ast} [\overleftarrow{1}_{\ast
i_1}|\widetilde{\alpha}_{i_1i_2}|\ldots|
\widetilde{\alpha}_{i_li_0}]})$ vanishes for $l$ large enough since
$\alpha$ is upper-triangular.

\medskip

Now Theorem \ref{eul} follows from the above three lemmas and the
following proposition:
\medskip
\begin{proposition}
Let $N_0$ be a free $A$-module. Then the map
$$\str:\rC_\bullet(\End_{\per A}(N_0))\to\rC_\bullet(A)$$ is a
quasi-isomorphism of complexes. Moreover, for any
$x\in\rHH_\bullet(\End_{\per A}(N_0))$ one has $x=\str(x)$ in
$\rHH_\bullet(\per A)$.
\end{proposition}
The proof of this statement is very similar to the proof of Theorem
1.2.4 from \cite{L} and we will omit it.

\medskip
\subsection{Computing the integral}\label{ci}
In Section \ref{hrr} we introduced an ``integral" $$\int:
\rHH_\bullet(\End_k(X))\to\rHH_\bullet(\per k)\simeq k$$ for any
complex of vector spaces $X$ with finite dimensional total
cohomology. In this section we will present an explicit formula for
this integral based on the results of \cite{FLS} (see also
\cite{R}). This, together with (\ref{wedge}), will give us an
explicit formula for computing the pairing (\ref{pairing}).

To exclude the trivial case, we will assume that $X$ has non-zero
cohomology.

Let us fix a pair of degree 0 maps $p:X\rightarrow\rH^\bullet(X)$
and $i:\rH^\bullet(X)\rightarrow X$ that establish the homotopy
equivalence between the complex $X$ and its cohomology
$\rH^\bullet(X)$:
$$
pi=1_{\rH^\bullet(X)},\quad ip=1_X-[d_X,H]
$$
where $H: X\to X$ is a degree $-1$ map.

Here is an explicit formula for the integral:
\medskip
\begin{theorem}\label{integral}
The following map is a quasi-isomorphism:
$$
\phi: \rC_\bullet(\End_k(X))\to k, \quad
T_1[T_2|\ldots|T_n]\mapsto\sum_{j=0}^{n-1}\str_{\rH^\bullet(X)}
(\cF_{n}(\tau^j(T_1[T_2|\ldots|T_n]))),
$$
where $\str_{\rH^\bullet(X)}$ is the ordinary super-trace,
$$
\tau(T_1[T_2|\ldots|T_n])=(-1)^{|sT_n|(|sT_1|+\ldots+|sT_{n-1}|)}T_n[T_1|\ldots|T_{n-1}],
$$
and $\cF_{n}: \End_k(X)^{\otimes n}\to \End_k(\rH^\bullet(X))$ is
given by
$$
\cF_{n}(T_1[T_2|\ldots|T_n])=
pT_1HT_2H\cdot\ldots\cdot HT_ni.
$$
Furthermore, the induced isomorphism $\rHH_\bullet(\End_k(X))\simeq
k$ coincides with $\int$.
\end{theorem}
Let us sketch the idea of the proof. That $\phi$ is a morphism of
complexes can be verified by a direct computation. Alternatively,
this follows from Lemma 2.4 of \cite{FLS} and the fact that the
collection $\cF_n$, $n=1,2,\ldots$, gives rise to an
$A_\infty$-morphism from the DG algebra $\End_k(X)$ to the DG
algebra $\End_k(\rH^\bullet(X))$. Moreover, the latter morphism is
an $A_\infty$-quasi-isomorphism, therefore
$\rHH_\bullet(\End_k(X))\simeq\rHH_\bullet(\End_k(\rH^\bullet(X)))\simeq
k$ which proves that $\phi$ is a quasi-isomorphism.

It remains to prove that the induced map $\rHH_\bullet(\End_k(X))\to
k$ coincides with $\int$. Obviously, it suffices to fix a non-zero
generator of $\rHH_0(\End_k(X))$ and to show that the values of both
functionals on this generator coincide. Let us start by describing a
generator of $\rHH_0(\End_k(X))$.

The endomorphism $ip$ is an idempotent. Let us denote its image by
$\Harm^\bullet(X)$. Clearly, $\Harm^\bullet(X)$ is a finite
dimensional subspace of $X$ isomorphic to $\rH^\bullet(X)$. Fix $n$
such that the component $\Harm^n(X)$ is non-zero and let $\pi$ stand
for the projection in $\Harm(X)$ onto this component parallel to
other graded components. Then the endomorphism $\Pi=\pi
ip\in\End^0_k(X)$ represents a non-zero element of
$\rHH_0(\End_k(X))$. It is immediate that
$\phi(\Pi)=(-1)^n\dim\,\rH^n(X)$.

On the other hand, $\Pi$ and $p\pi i$ define the same element of
$\rHH_0(\per k)$ ($p\pi i$ is just for the projection in
$\rH^\bullet(X)$ onto the component $\rH^n(X)$ parallel to other
graded components). Indeed, $\Pi-p\pi i=\pi ip-p\pi i=b(\pi i[p])$.
To finish the proof, observe that the element of $\rHH_0(\per k)$,
defined by $p\pi i$, coincides with the one, defined by
$(-1)^n\dim\,\rH^n(X)\cdot 1\in\End_k(k)$.

\medskip

Let us point out a couple of straightforward corollaries of Theorem
\ref{integral} and formula (\ref{pairing}).

Let $A$ be a finite dimensional associative algebra. Then its
Hochschild homology groups $\rHH_\bullet(A)$ are concentrated in
non-positive degrees. Therefore, among the pairings $\langle
\,,\,\rangle: \rHH_n(A)\times\rHH_{-n}(A^\op)\to k$, only the one
corresponding to $n=0$ survives. In this case we have

\medskip
\begin{corollary}\label{assocpair}
For an associative algebra $A$, the pairing $$\langle \,,\,\rangle:
A/[A,A]\times A^\op/[A^\op,A^\op]\to k$$ is given by
\begin{eqnarray*}\label{assocpairing}
\langle\, a, b\,\rangle=\tr_A(L(a)R(b)).
\end{eqnarray*}
(In the right-hand side, $a$ and $b$ stand for elements of $A$ and
$A^\op$, respectively, and in the left-hand side  $a,b$ stand for
the corresponding classes in the Hochschild homology.)
\end{corollary}
\medskip

Let now $A$ be a finite dimensional graded algebra. Since $A$ is
equipped with the trivial differential, we can set $H=0$ in Theorem
\ref{integral} and obtain

\medskip
\begin{corollary}\label{gradpair}
For a graded $A$, the pairing of two cycles
$$
a=a_0+\sum a'_0[a'_1]+\sum a''_0[a''_1|a''_2]+
\ldots\in\rC_\bullet(A),
$$
$$
b=b_0+\sum b'_0[b'_1]+\sum
b''_0[b''_1|b''_2]+\ldots\in\rC_\bullet(A)
$$
is given by
$$
\langle\,a,\,b\,\rangle=\str_{A}(L(a_0)R(b_0)).
$$
\end{corollary}

\medskip
\section{Examples}\label{eg}
\subsection{Directed algebras}\label{dir}

In this section, we describe how the Hirzebruch-Riemann-Roch formula
looks like for a special class of finite dimensional associative
algebras.

Let $\cV$ be a $k$-linear category with finite number of objects,
say $\{v_s\}_{s\in S}$, and finite dimensional $\Hom$-spaces.
Suppose there is a bijection
$$f:\{1,2,\ldots,n\}\to I$$ such that
\begin{equation}\label{a}
\Hom_{\cV}(v_{f(i)},v_{f(j)})=\begin{cases}k & i=j\\ 0 & i>j \end{cases}.
\end{equation}
Of course, $f$ doesn't have to be unique. Let us denote the algebra
of this category by $A(\cV)$:
$$
A(\cV)=\bigoplus_{s,t\in S}\Hom_{\cV}(v_s,v_t).
$$
We will call such algebras (as well as the underlying categories)
{\it directed}.

Let us denote the abelian category of finite dimensional right
$A(\cV)$-modules by $\mmod A(\cV)$. The following simple result is
very well known.

\medskip
\begin{proposition}\label{resol}
Any module $N\in\mmod A(\cV)$ admits a projective resolution of
finite length.
\end{proposition}
Let us prove this. Fix a map $f$ as above and denote $1_{v_{f(i)}}$
simply by $1_i$. Denote also the projective modules $1_i\,A(\cV)$ by
$P_i$. Clearly,
$$\dim\Hom_{\mmod
A(\cV)}(P_i,P_j)=\dim\Hom_{\cV}(v_{f(i)},v_{f(j)}).$$  Thus, by
(\ref{a})
\begin{equation}\label{b}
\dim\Hom_{\mmod A(\cV)}(P_i,P_j)=\begin{cases}k & i=j\\ 0 & i>j \end{cases}.
\end{equation}
Fix $N\in\mmod A(\cV)$. The canonical morphism
$$
p: \bigoplus_{i=1}^n \Hom_{\mmod A(\cV)}(P_i,N)\otimes_kP_i\to N
$$
is surjective. The kernel of this morphism satisfies the property
$$\Hom_{\mmod A(\cV)}(P_n,\mathsf{Ker}\, p)=0.$$ To see this, apply the functor
$\Hom_{\mmod A(\cV)}(P_n,\, -\,)$ to the short exact sequence
$$
0\to \mathsf{Ker}\, p\to \bigoplus_{i=1}^n \Hom_{\mmod A(\cV)}(P_i,N)\otimes_kP_i\to N \to0
$$
and use the property (\ref{b}).

To finish the proof, apply the same argument to $\mathsf{Ker}\, p$
instead of $N$ etc.

\medskip

Observe that $\rHH_0(A(\cV))$ is spanned by the idempotents
$1_{v_s}$, $s\in S$ (or rather their classes in the quotient
$A(\cV)/[A(\cV),A(\cV)]$). In terms of these elements, the pairing
$\langle\,,\,\rangle$ on $\rHH_0(A(\cV))\times\rHH_0(A(\cV)^\op)$ is
given by
$$
\langle\,1_t\,,\,1_s^\vee\,\rangle=\dim\Hom_{\cV}(v_s,v_t).
$$

Let us derive the Hirzebruch-Riemann-Roch formula for finite
dimensional modules over directed algebras. It is well known and was
obtained in \cite[Section 2.4]{Ri}.

Let us keep the notations from the proof of Proposition
\ref{direct}. Set $$d_{ij}:=\dim\Hom_{\cV}(v_{f(i)},v_{f(j)}).$$ Let
$M,N\in\mmod A(\cV)$. As we saw above, $M$ and $N$ admit finite
length resolutions by direct sums of the projective modules $P_i$.
Let us fix two such resolutions $P(M)$ and $P(N)$. We know that
$\eul(P(M))$, $\eul(P(N))$ are linear combinations of $1_i$'s:
$$
\eul(P(M))=\sum_{i=1}^n a_i\cdot 1_i, \quad \eul(P(N))=\sum_{i=1}^n b_i\cdot 1_i.
$$
Since $1_j=\eul(P_j)$, we have
\begin{eqnarray*}
(\underline{\dim} M)_j:=\Hom_{\mmod A(\cV)}(P_j,M)=\Hom_{\Ho(\Mod
A(\cV))}(P_j,P(M))=\langle\,\eul(P(M))\,,\,1_j^\vee\,\rangle \\
=\sum_{i=1}^n d_{ji}a_i
\end{eqnarray*}
and similarly $(\underline{\dim} N)_j=\sum_{i=1}^n d_{ji}b_i$.
Therefore,
\begin{eqnarray*}
\sum_{l}(-1)^l\dim\mathrm{Ext}_{\mmod A(\cV)}^l(M,N)=\chi(P(M),P(N))=
\langle\,\eul(P(N))\,,\,\eul(P(M))^\vee\,\rangle\\
=\sum_{i,j}b_ia_jd_{ji}.
\end{eqnarray*}
Since $ a_j=\sum_k(d^{-1})_{jk}(\underline{\dim} M)_k, \quad
b_i=\sum_k(d^{-1})_{il} (\underline{\dim} N)_l, $ we get the
following generalization of Ringel's formula:
\begin{equation}\label{grf}
\sum_{l}(-1)^l\dim\mathrm{Ext}_{\mmod A(\cV)}^l(M,N)=\sum_{i,j}(\underline{\dim} M)_i(d^{-1})_{ij}
(\underline{\dim} N)_j.
\end{equation}

\medskip
\subsection{Proper noncommutative DG-schemes arising from orbifold singularities}\label{orbi}

In this section, we will describe certain proper DG
algebras\footnote{All of them are DG algebras with the trivial
differential.} which arise from quotient singularities of the form
$\C^n/G$, where $G$ is a finite group.

Let $V=\C^n$ be a finite dimensional complex vector space and $G$ a
finite subgroup of $SL(V)\cong SL_n(\C)$. Then $G$ acts on the
polynomial algebra $\C[V]$ via $(gf)(x)=f(g^{-1}x)$. The spectrum
$X=V/G$ of the algebra $\C[V]^G$ of $G$-invariant polynomials is a
singular affine variety. The central problem in the study of such
singular varieties is to construct their ``most economical"
resolutions, which are called {\it crepant}: a resolution $\pi: Y\to
X$ is crepant, if $\pi$ preserves the canonical classes\footnote{A
crepant resolutions of $X$, if exists, is a noncompact Calabi-Yau
variety since the top-degree form on $V$ is $G$-invariant and
therefore the canonical sheaves of $X$ and $Y$ are trivial.}, i.e.
$\pi^*(\omega_X)=\omega_Y$.

The derived Mckay correspondence \cite{R1,R2,KV,BKR} is a program
around the following conjecture and various versions thereof:

\medskip
\begin{quote}
\emph{For any crepant resolution $Y\to X$, the bounded derived
category $D(Y)$ of coherent sheaves on $Y$ is equivalent to the
bounded derived category $D^G(V)$ of $G$-equivariant coherent
sheaves on $V$.}
\end{quote}
\medskip

\noindent In other words, all crepant resolutions of a fixed
singularity are expected to be isomorphic as noncommutative
DG-schemes. The conjecture is known to be true for finite subgroups
of $SL(2)$ \cite{KV} and $SL(3)$ \cite{BKR} (see also \cite{BKa} for
a result in higher dimensions).

Denote by $D^G_0(V)$ the subcategory in $D^G(V)$ of complexes
supported at the origin $0\in V$ and by $D_0(Y)$ the subcategory in
$D(Y)$ of complexes supported at the exceptional fiber $\pi^{-1}(0)$
(in the latter formula $0$ stands for the image of the origin of $V$
under the canonical projection $V\to X$). Then the above equivalence
of categories should induce an equivalence between $D_0(Y)$ and
$D^G_0(V)$ \cite{BKR}.

The Ext groups between any two objects of $D^G_0(V)$ vanish in all
but finitely many degrees and, thus, we are dealing with a proper
noncommutative DG-scheme. This scheme is the main subject of the
section.

Following \cite[Section 6.2]{HK}, consider the cross-product
$\Lambda(V,G)$ of the exterior algebra $\Lambda V$ and the group
algebra of $G$. In other words, as a vector space $\Lambda(V,G)$ is
the tensor product $\Lambda V\otimes \C[G]$. The product of two
elements is given by
$$
(v\otimes g)(w\otimes h) = (v\wedge g(w))\otimes gh, \quad v,w\in \Lambda V,\, g,h\in G.
$$
Equip $\Lambda(V,G)$ with the unique grading such that $\deg v=1$
and $\deg g=0$ for any $v\in V$ and $g\in G$. We will view
$\Lambda(V,G)$ as a DG algebra with the trivial differential.

The following conjecture is motivated by \cite{HK}:

\medskip
\noindent{\bf Conjecture.} {\it There is an equivalence of
triangulated categories $$D^G_0(V)\cong\Ho(\per \Lambda(V,G)).$$}
\medskip

Here is how the conjecture might be proved. The category $D^G_0(V)$
seems to be equivalent to the category $D^b(f.d.\,\C[V]\rtimes G)$,
where $\C[V]\rtimes G$ is the cross-product of the polynomial
algebra and the group algebra of $G$ and $f.d.\,\C[V]\rtimes G$ is
the abelian category of finite dimensional graded $\C[V]\rtimes
G$-modules. Every such module admits a finite filtration by simple
$\C[V]\rtimes G$-modules. The latter are the $\C[V]\rtimes
G$-modules obtained from simple $\C[G]$-modules via ``restriction of
scalars"
$$
\C[V]\rtimes G\to\C[G],\quad f(x)\otimes g\mapsto f(0)g,\quad f(x)\in\C[V],\, g\in G.
$$
Let us denote the simple $\C[V]\rtimes G$-module, corresponding to
an irreducible representations $\rho$ of $G$, by $S_\rho$. Then,
using the technique described in \cite{K4}, we may conclude that
$D^b(f.d.\,\C[V]\rtimes G)$ is equivalent to the category $\Ho(\per
\cA)$ for some $A_\infty$ algebra $\cA$ with
$$
\rH^\bullet(\cA)=\Ext^\bullet(\oplus_\rho S_\rho,\oplus_\rho S_\rho),
$$
where the sum in the right-hand side is taken over irreducible
representations of $G$. According to \cite[Section 6.2]{HK}, the
algebra $\C[V]\rtimes G$ is quadratic and Koszul, and its Koszul
dual is exactly $\Lambda(V,G)$. Then, by \cite[Section 2.2]{K4}, the
$A_\infty$ algebra $\cA$ is formal. Finally, we expect that
$\Ext^\bullet(\oplus_\rho S_\rho,\oplus_\rho S_\rho)$ is Morita
equivalent to $\Lambda(V,G)$.

\medskip

Whether the conjecture is true or not, it is clear that the
algebraic triangulated categories of the form $\Ho(\per
\Lambda(V,G))$ should play a role in the study of the quotient
singularities.

\medskip

Let us compute the pairing $\langle\,,\,\rangle$ on
$\rHH_0(\Lambda(V,G))\times\rHH_0(\Lambda(V,G)^\op)$.

We start by noticing that, in general, the space
$\rHH_0(\Lambda(V,G))$ is infinite dimensional (this is already so
in the simplest case $V=\C$, $G=\{1\}$). However, the pairing
$\langle\,,\,\rangle$ vanishes on a subspace of finite codimension
(this follows from Corollary \ref{gradpair}). In fact, the pairing
is determined by its restriction onto the finite dimensional
subspace
$$\rHH_0(\C[G])\times \rHH_0(\C[G]^\op)\subset\rHH_0(\Lambda(V,G))\times \rHH_0(\Lambda(V,G)^\op).$$
(Here we are using the natural embedding $\C[G]\to \Lambda(V,G)$
which induces an embedding $\rHH_0(\C[G])\to \rHH_0(\Lambda(V,G))$.)
Furthermore, it is well known that $\rHH_0(\C[G])$ is spanned by
(the homology classes of) the characters of irreducible
representations of $G$. Let us denote the character of an
irreducible representation $\rho$ by $\chi_\rho$:
$$
\chi_\rho=\sum_g \tr(\rho(g))g.
$$
Using basic harmonic analysis on $G$, it is easy to show that the
element $\pi_\rho=\frac{\dim \rho}{|G|}\chi_\rho$ is an idempotent
in $\Lambda(V,G)$ (it is nothing but the Euler class of the DG
$\Lambda(V,G)$-module $\pi_\rho\cdot\Lambda(V,G)$). Thus, we have to
compute $$\langle\,\pi_{\rho_1},\pi^\vee_{\rho_2}\,\rangle=
\str_{\Lambda(V,G)}(L(\pi_{\rho_1})R(\pi_{\rho_2}))$$ for two
irreducible representations $\rho_1$, $\rho_2$.

Let $W$ be the space of some representation of $G$. Then
$W\otimes\C[G]$ carries a natural $\C[G]$-bimodule structure,
defined as follows:
$$
g(w\otimes h)k = g(w)\otimes ghk, \quad w\in W,\, g,h,k\in G.
$$
In particular, the graded components
$\Lambda^n(V,G)=\Lambda^nV\otimes\C[G]$ of the algebra
$\Lambda(V,G)$ are $\C[G]$-bimodules and we have
\begin{eqnarray*}
\str_{\Lambda(V,G)}(L(\pi_{\rho_1})R(\pi_{\rho_2}))=\sum_{n=0}^{\dim
V}(-1)^n\,\tr_{\Lambda^n(V,G)}(L(\pi_{\rho_1})R(\pi_{\rho_2}))\\=\sum_{n=0}^{\dim
V}(-1)^n\,\dim(\pi_{\rho_1}\Lambda^n(V,G)\pi_{\rho_2}).
\end{eqnarray*}
Therefore, we will start by computing
$\dim(\pi_{\rho_1}(W\otimes\C[G])\pi_{\rho_2})$ for an arbitrary
$W$.

Let us introduce a matrix $d^W$ of non-negative integers by the
following formula:
$$
W\otimes\rho=\bigoplus_\sigma d^W_{\sigma\rho}\,\sigma,
$$
where $\rho$ and $\sigma$ run through the set of irreducible
representations of $G$. Let us denote the representation, dual to
$\rho$, by $\rho'$. Then, as a $\C[G]$-bimodule
$$
W\otimes\C[G]=\bigoplus_\rho(W\otimes
\rho)\boxtimes\rho'=\bigoplus_{\rho,\sigma}d^W_{\sigma\rho}\,\sigma\boxtimes\rho'.
$$
Thus,
$$
\dim(\pi_{\rho_1}(W\otimes\C[G])\pi_{\rho_2})=\dim\rho_1\,\dim\rho_2\,d^W_{\rho_1\rho_2},
$$
which gives us the following formula for
$\langle\,\pi_{\rho_1},\pi^\vee_{\rho_2}\,\rangle$:
\begin{equation}\label{orbhrr}
\langle\,\pi_{\rho_1},\pi^\vee_{\rho_2}\,\rangle=\dim\rho_1\,\dim\rho_2\,\sum_{n=0}^{\dim
V}(-1)^n\,d^{\Lambda^nV}_{\rho_1\rho_2}.
\end{equation}

\medskip
\section{More on the pairing $\langle\,,\,\rangle$}\label{pairi}
\subsection{Smooth proper DG algebras}\label{smooth}

Recall \cite{KS} that a DG algebra is said to be (homologically)
smooth if there is a perfect right DG $A^\op\otimes A$-module $P(A)$
together with a quasi-isomorphism $P(A)\to A$ of right DG
$A^\op\otimes A$-modules.

To have an example at hand, observe that
\medskip
\begin{proposition}\label{direct}
Any directed algebra is smooth.
\end{proposition}
Indeed, it is clear that $A(\cV)^\op\otimes A(\cV)\cong
A(\cV^\op\otimes\cV)$. Therefore, by Proposition \ref{resol}, any
finite dimensional $A(\cV)^\op\otimes A(\cV)$-module admits a finite
projective resolution. What remains is to apply this to $A(\cV)$ and
observe that any finite complex of projective bimodules over an
associative algebra is a perfect DG bimodule in our sense.

\medskip

The aim of this section is to prove that the pairing
$$\langle\,,\,\rangle: \rHH_n(\per A)\times\rHH_{-n}(\per A^\op)\to
k,$$ is non-degenerate for any proper smooth DG algebra $A$.
The proof is based on the observation that the pairing is inverse to
the Euler class $\eu(A)$ of the $A$-bimodule $A$. The author learned
about this idea from \cite{KS1}\footnote{Although \cite{KS1} is
still in preparation, the argument with the Euler class of the
``diagonal" is already well known among the experts \cite{Ka}.}.

\medskip
\begin{theorem}\label{t}
Let $A$ be a proper smooth DG algebra. Then the pairing
$\langle\,,\,\rangle$ is non-degenerate.
\end{theorem}
Indeed, fix a perfect resolution $P(A)\stackrel{p}\to A$ in the
category of right DG $A^\op\otimes A$-modules. Then, for any right
perfect DG $A$-module $X$, we have a morphism
$$
1\otimes p:X\otimes_AP(A)\to X\otimes_AA\simeq X.
$$
By Proposition \ref{flat}, $1\otimes p$ is a quasi-isomorphism. On
the other hand, by Proposition \ref{tensor}, both $X\otimes_AP(A)$
and $X$ are perfect and, in particular, homotopically projective. It
is well known that a quasi-isomorphism between two homotopically
projective modules is actually a homotopy equivalence (see, for
instance, the proof of Lemma 10.12.2.2 in \cite{BL}). Thus,
$1\otimes p:X\otimes_AP(A)\to X\otimes_AA\simeq X$ is a homotopy
equivalence.

What we have just proved is that the quasi-isomorphism
$P(A)\stackrel{p}\to A$ gives rise to a weak homotopy equivalence of
the DG functors $T_{P(A)}\to\mathrm{Id}_{\per A}$ where
$\mathrm{Id}$ stands for the identity endofunctor. Then, as a
corollary of Theorem \ref{hominv}, we get the following result: the
linear map $\rHH(T_{P(A)}): \rHH_\bullet(\per A)\to
\rHH_\bullet(\per A)$ coincides with the identity map. On the other
hand, by Theorem \ref{main}, the map $\rHH(T_{P(A)})$ is given by
the 'convolution' with $\eu'(P(A))$, so the convolution with
$\eu'(P(A))$ is the identity map. This proves that the left kernel
of the pairing is trivial, i.e. for any $n$ we have an embedding
$$
\rHH_n(\per A)\to\rHH_{-n}(\per A^\op)^*.
$$
One of the results of \cite{Sh} says that the Hochschild homology of
an arbitrary proper smooth DG algebra is finite dimensional. Thus,
to prove that the right kernel of the pairing is trivial, it is
enough to show that $\dim\rHH_n(\per A)=\dim\rHH_{-n}(\per A^\op)$.
This can be done by replacing $A$ by $A^\op$ in the above argument.

\medskip

Let us point out one interesting corollary of this
result\footnote{If $k$ is perfect, this result also follows from
Proposition 2.5 of \cite{Ke-1} and Morita invariance of the
Hochschild homology.}:
\medskip
\begin{corollary}\label{c}
If $A$ is a smooth proper associative algebra then
$$
\rHH_n(A)=\begin{cases}A/[A,A] & n=0\\0 & \text{otherwise}\end{cases}
$$
\end{corollary}
Indeed, the Hochschild homology of such an algebra is concentrated
in non-positive degrees. Thus, by the non-degeneracy of the pairing,
the Hochschild homology groups, sitting in negative degrees, have to
vanish.

\medskip
This corollary, together with Proposition \ref{direct}, implies
$\rHH_n(A(\cV))=0$ for any directed algebra $A(\cV)$ and any
$n\neq0$. This result was obtained by a different method in
\cite{Ci}.

Another application of the corollary is related to the so-called
noncommutative Hodge-to-de Rham degeneration conjecture. Roughly
speaking, the conjecture claims that the $B$-operator
$B:\rHH_\bullet(A)\to\rHH_{\bullet-1}(A)$ (see \cite{GJ,Ts})
vanishes whenever $A$ is proper and smooth. It was formulated, in a
stronger form, by M. Kontsevich and Y. Soibelman \cite{KS} and
proved, in the partial case of DG algebras concentrated in
non-negative degrees, by D. Kaledin \cite{Ka}. The above corollary
implies the conjecture in the case of algebras with the trivial
differential and grading.

\medskip
\subsection{Relation to Topological Field Theories}\label{TCFT}
This section is devoted to yet another application of our results.
Namely, we will discuss the relevance of the pairing (\ref{pairing})
to the Topological Field Theories (TFT's) constructed in \cite{Cos,
KS}.

To begin with, let us recall that by a {\it trace} on a DG algebra
$A$ one understands a (homogeneous) functional $\tau: A\to k$ such
that
$$
\tau(da)=0,\quad \tau([a,b])=0,\quad a,b\in A.
$$

Let $A$ be a proper DG algebra. Suppose the algebra possesses a
degree $-d$ trace $\tau$ satisfying the following condition: the
induced degree $-d$ pairing
$$
\rH^\bullet(A)\times\rH^\bullet(A)\to k,\quad (a,b)\mapsto \tau(ab)
$$
is non-degenerate. Then the pair $(A,\tau)$ is called a
$d$-dimensional (compact) Calabi-Yau DG algebra
\cite{KS}\footnote{Actually, the authors of \cite{Cos, KS} work with
CY $A_\infty$ algebras and categories.}. Sometimes, we will write
$A$ instead of $(A,\tau)$.

Observe that the algebra $\Lambda(V,G)$, we studied in Section
\ref{orbi}, carries a natural structure of a $\dim\,V$-dimensional
CY DG algebra. Namely, fix a non-zero element
$\omega\in\Lambda^{\dim\,V}V$ and set \cite{HK}:
$$
\tau_{\omega}(v\otimes g)=\begin{cases} 0 & v\in \Lambda^nV,\, n<\dim\,V\\ \delta_{1g} &
v=\omega\end{cases}.
$$

Before we proceed any further, we would like to mention that there
is a different class of CY algebras whose theory is now being
actively developed \cite{Gi1}. These latter CY algebras are
noncommutative analogs of {\it noncompact} smooth CY varieties (a
good example of such an algebra is the cross-product $\C[V]\rtimes
G$ we mentioned in Section \ref{orbi}).

A $d$-dimensional TFT is defined as follows (we refer to \cite{Cos,
KS} for details). Let $\cM(n,m)=\bigcup_{g\geq0}\cM_g(n,m)$ denote
the moduli space of Riemann surfaces with $n$ incoming and $m$
outgoing boundaries ($g$ denotes the genus). Fix a graded vector
space $\bH_\bullet$. By definition, $\bH_\bullet$ carries a
structure of a $d$-dimensional TFT if one has a collection of linear
maps
\begin{equation}\label{action}
\rH_\bullet(\cM(n,m))\otimes\bH_\bullet^{\otimes
n}\to\bH_\bullet^{\otimes m},\quad n\geq1, \,\,m\geq0
\end{equation}
(here $\rH_\bullet$ in the left-hand side denotes the singular
homology) satisfying the following conditions:

\noindent (1) the maps are compatible with the operation
$$
\cM(m,l)\times\cM(n,m)\to \cM(n,l)
$$
of gluing of two surfaces along the boundary components and the
operation
$$
\cM(n,m)\times \cM(p,q)\to \cM(n+p,m+q)
$$
of taking the disjoint union of surfaces;

\noindent (2) elements of $\rH_\bullet(\cM_g(n,m))$ act by operators
of degree $d(2-2g-n-m)$.

One has \cite{Cos, KS}:

\medskip

\begin{quote}
\emph{For any $d$-dimensional CY DG algebra $A$ the Hochschild
homology $\rHH_\bullet(A)$ carries a canonical structure of
$d$-dimensional TFT.}\footnote{In fact, a much stronger result is
obtained in \cite{Cos, KS}, namely, that the action (\ref{action})
exists on the level of complexes that compute the singular homology
of the moduli spaces and the Hochschild homology of the algebra. }
\end{quote}

\medskip

One immediate consequence of this result is that there is a natural
degree 0 pairing -- let us denote it by
$\langle\,\,,\,\,\rangle_\tau$ -- on the Hochschild homology of a
$d$-dimensional CY DG algebra, given by a generator of
$\rH_0(\cM_0(2,0))$. The following conjecture relates this pairing
to the one constructed in the present work\footnote{This conjecture
was suggested to the author by Y. Soibelman and K. Costello.}:

\medskip
\noindent{\bf Conjecture.} {\it For any CY DG algebra $A$, the
pairing $\langle\,\,,\,\,\rangle_\tau$ coincides with the pairing
(\ref{pairing}), i.e. for any $a,b\in\rHH_\bullet(A)$
\begin{equation}\label{conjec}
\langle\,a,b\,\rangle_\tau=\langle\,a,b^\vee\,\rangle,
\end{equation}
where $\,^\vee$ is the isomorphism
$\rHH_\bullet(A)\to\rHH_\bullet(A^\op)$ defined by (\ref{aaop}).}

\medskip
We note that this conjecture, together with Theorem \ref{t}, would
imply the following result conjectured in \cite[Section 11.6]{KS}:

\medskip
\noindent{\bf Corollary.} {\it For any smooth CY DG algebra $A$, the
pairing $\langle\,\,,\,\,\rangle_\tau$ is non-degenerate.}
\medskip

To present a piece of evidence in favor of the conjecture, let us
prove it in the case of an associative Calabi-Yau algebra, when the
grading and the differential are both trivial. Observe that such a
CY algebra is nothing but a symmetric Frobenius algebra \cite{Ko}.

To compute the left-hand side of (\ref{conjec}), we will use an
explicit description of the action (\ref{action}) based on graphs
\cite[Section 11.6]{KS}. In the language of \cite{KS}, the pairing
$\langle\,,\,\rangle_\tau$ corresponds to the following graph:

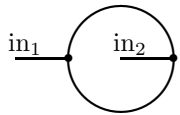
\setlength{\unitlength}{1cm}
\begin{picture}(2,2)
\thicklines \put(7.5,1){\circle{1.4}} \thicklines
\put(7.5,1){\line(1,0){0.7}} \thicklines
\put(6.1,1){\line(1,0){0.7}} \put(6.8,1){\circle*{0.1}}
\put(8.2,1){\circle*{0.1}} \put(6.0,1.1){\small{in}$_1$}
\put(7.4,1.1){\small{in}$_2$}
\end{picture}

Let us fix a symmetric Frobenius algebra $A=(A,\tau)$. Since $A$ is
finite dimensional and the bilinear form $\tau(ab)$ is
non-degenerate, there exists a unique element
$$
\Phi=\sum_k\phi'_k\otimes \phi''_k\otimes \phi'''_k\in A\otimes A\otimes A
$$
satisfying the property
\begin{equation}\label{phi}
\tau(abc)=\sum_k \tau(a\phi'_k)\tau(b\phi''_k)\tau(c\phi'''_k)
\end{equation}
for every $a,b,c\in A$. Notice that $\Phi$ is cyclically symmetric
because $\tau(ab)$ is symmetric. According to \cite{KS},
$\langle\,a,b\,\rangle_\tau$ can be computed by means of the above
graph as follows: attach $a$ to the vertex marked in$_1$ and $b$ to
the vertex marked in$_2$; attach two copies of the tensor $\Phi$ to
the remaining two vertices; contruct all the tensors along all four
edges of the graph, using the pairing $a\times b\mapsto\tau(ab)$.
Here is the result:
\begin{equation*}\label{formulaforpair}
\langle\,a,b\,\rangle_\tau=\sum_{k,l}\tau(a\phi'_k)\tau(b\phi'_l)\tau(\phi''_k\phi''_l)
\tau(\phi'''_k\phi'''_l).
\end{equation*}

By (\ref{phi}), the latter formula can be simplified as follows:
$$
\langle\,a,b\,\rangle_\tau=\sum_{k}\tau(a\phi'_k)\tau(b\phi''_k\phi'''_k).
$$
To simplify the formula further, consider the unique symmetric
element $$\gamma=\sum_i \gamma'_i \otimes \gamma''_i\in A\otimes A$$
satisfying the property
\begin{equation}\label{symtens}
a=\sum_i \gamma'_i\tau(\gamma''_ia)=\sum_i\tau(a\gamma'_i)\gamma''_i
\end{equation}
for every $a\in A$. Then it is easy to see that $ \Phi=\sum_{i,j}
\gamma'_i\otimes \gamma'_j\gamma''_i\otimes \gamma''_j. $ Indeed,
\begin{eqnarray*}
\sum_{i,j}\tau(a\gamma'_i)\tau(b\gamma'_j\gamma''_i)\tau(c\gamma''_j)=\sum_i\tau(a\gamma'_i)\tau(b\sum_j\gamma'_j\tau(c\gamma''_j)\gamma''_i)
=\sum_i\tau(a\gamma'_i)\tau(bc\gamma''_i)\\
=\tau(a\sum_i\gamma'_i\tau(bc\gamma''_i))=\tau(abc).
\end{eqnarray*}
Thus,
\begin{eqnarray*}
\langle\,a,b\,\rangle_\tau=\sum_{k}\tau(a\phi'_k)\tau(b\phi''_k\phi'''_k)=\sum_{i,j}\tau(a\gamma'_i)\tau(b\gamma'_j\gamma''_i\gamma''_j)
=\sum_{j}\tau(b\gamma'_j\sum_i\tau(a\gamma'_i)\gamma''_i\gamma''_j)\\=\sum_{j}
\tau(b\gamma'_ja\gamma''_j).
\end{eqnarray*}
Since $\gamma$ is symmetric, we arrive at the following formula
$$
\langle\,a,b\,\rangle_\tau=\sum_i\tau(a\gamma'_ib\gamma''_i).
$$

By Corollary \ref{assocpair}, we have
$\langle\,a,b^\vee\,\rangle=\tr_A(L(a)R(b))$. Thus, for a symmetric
Frobenius algebra, the above conjecture boils down to the following
identity:
$$
\sum_i\tau(a\gamma'_ib\gamma''_i)=\tr_A(L(a)R(b)), \quad a,b\in A.
$$
To prove it, we observe that under the canonical isomorphism
$\End_k(A)\cong A\otimes A^*$ the operators $L(a)$, $R(b)$ get
mapped to the elements
$$
\sum_i a\gamma'_i\otimes \tau(\gamma''_i\cdot-), \quad \sum_j \gamma'_jb\otimes
\tau(\gamma''_j\cdot-),
$$
respectively (this follows from the definition (\ref{symtens}) of
$\gamma$). Therefore,
\begin{eqnarray*}
\tr_A(L(a)R(b))=\sum_{i,j}\tau(\gamma''_ja\gamma'_i)\tau(\gamma''_i\gamma'_jb)=\sum_i\tau(\gamma''_i\gamma'_j\sum_j\tau(\gamma''_ja\gamma'_i)b)\\=
\sum_i\tau(\gamma''_ia\gamma'_ib)=\sum_i\tau(a\gamma'_ib\gamma''_i),
\end{eqnarray*}
which finishes the proof.

\medskip
The same proof should work for graded CY DG algebras with the
trivial differential.

\medskip
\appendix
\section{Proof of Proposition \ref{veeiso}}\label{proof}
Clearly, the morphism (\ref{aaop}) is invertible. We have to show
that it commutes with the differentials. It is obvious that
$\,^\vee$ respects the first differential $b_0$ as its definition
doesn't involve multiplication. Let us show by a direct computation
that $\,^\vee$ commutes with the second differential $b_1$. Let us
denote the multiplication in $A^\op$ by $\ast$. To simplify
computations, we will also use the notations
$\xi_i=|a_0|+|sa_n|+|sa_{n-1}|+\ldots+|sa_{i+1}|$ and
$f(a_1,a_2,\ldots,a_n)=\sum_{1\leq i<j\leq n}|sa_i||sa_j|$. One has:
\begin{eqnarray*}
b_1((a_0[a_1|a_2|\ldots |a_n])^\vee)=(-1)^{n+f(a_1,a_2,\ldots,a_n)}b_1(a_0[a_n|a_{n-1}|
\ldots |a_1])\\
=(-1)^{n+f(a_1,a_2,\ldots,a_n)}((-1)^{|a_0|}a_0\ast a_n[a_{n-1}|\ldots |a_1]\\ +\sum\limits_{i=1}^{n-1}(-1)^{\xi_{i}}a_0[a_n|a_{n-1}|\ldots |a_{i+1}\ast a_{i}|\ldots|a_1]\\
-(-1)^{\xi_{1}(|a_1|+1)}a_1\ast a_0[a_{n}|a_{n-1}|\ldots |a_2])\\
=(-1)^{n+f(a_1,a_2,\ldots,a_n)}((-1)^{|a_0|+|a_0||a_n|}a_na_0[a_{n-1}|\ldots |a_1]\\ +\sum\limits_{i=1}^{n-1}(-1)^{\xi_{i}+|a_{i+1}||a_{i}|}a_0[a_n|a_{n-1}|\ldots |a_{i}a_{i+1}|
\ldots|a_1]\\
-(-1)^{\xi_{1}(|a_1|+1)+|a_1||a_0|}a_0a_1[a_{n}|a_{n-1}|\ldots
|a_2])
\end{eqnarray*}
On the other hand,
\begin{eqnarray*}
(b_1(a_0[a_1|a_2|\ldots |a_n]))^\vee=(-1)^{|a_0|}(a_0a_1[a_2|\ldots |a_n])^\vee\\ +\sum\limits_{i=1}^{n-1}(-1)^{\eta_{i}}(a_0[a_1|a_2|\ldots |a_ia_{i+1}|\ldots|a_n])^\vee\\
-(-1)^{\eta_{n-1}(|a_n|+1)}(a_na_0[a_1|a_2|\ldots |a_{n-1}])^\vee\\
=(-1)^{|a_0|}(-1)^{n-1+f(a_2,\ldots,a_n)}a_0a_1[a_n|\ldots |a_2]\\ +\sum\limits_{i=1}^{n-1}(-1)^{\eta_{i}}(-1)^{n-1+f(a_1,a_2,\ldots,a_ia_{i+1},\ldots, a_n)}
a_0[a_n|a_{n-1}|\ldots |a_ia_{i+1}|\ldots|a_1]\\
-(-1)^{\eta_{n-1}(|a_n|+1)}(-1)^{n-1+f(a_1,a_2,\ldots,a_{n-1})}a_na_0[a_{n-1}|a_{n-2}|\ldots
|a_1]
\end{eqnarray*}
What remains is to compare the signs, i.e. to show that
$$
(-1)^{f(a_1,a_2,\ldots,a_n)}(-1)^{|a_0|+|a_0||a_n|}=(-1)^{\eta_{n-1}(|a_n|+1)}
(-1)^{f(a_1,a_2,\ldots,a_{n-1})},
$$
$$
(-1)^{f(a_1,a_2,\ldots,a_n)}(-1)^{\xi_{i}+|a_{i+1}||a_{i}|}=-(-1)^{\eta_{i}}
(-1)^{f(a_1,a_2,\ldots,a_ia_{i+1},\ldots,
a_n)},
$$
$$
(-1)^{f(a_1,a_2,\ldots,a_n)}(-1)^{\xi_{1}(|a_1|+1)+|a_1||a_0|}=(-1)^{|a_0|}
(-1)^{f(a_2,\ldots,a_n)},
$$
which is an easy computation.

\small{\it Mathematics Department, Kansas State University}

\small{\it 138 Cardwell Hall}

\small{\it Manhattan, KS 66506-2602}

\small{\it e-mail: shklyarov@math.ksu.edu}
\end{document}